\documentclass{ifacconf}

\usepackage{graphicx}     
\usepackage{natbib}        
\usepackage{mathtools}
\usepackage{amssymb, amsmath, latexsym}
\usepackage{algorithm}\usepackage{algpseudocode}
\usepackage{nicefrac}
\usepackage{color}

\usepackage[colorinlistoftodos,bordercolor=orange,backgroundcolor=orange!20,linecolor=orange,textsize=scriptsize]{todonotes}

\def \RR {\mathbb R}
\def \R {\mathbb R}
\def\eqdef{\overset{\text{def}}{=}}

\newcommand{\argmin}{\mathop{\arg\!\min}}
\newcommand{\argmax}{\mathop{\arg\!\max}}

\newcommand{\cN}{{\cal N}}

\newcommand{\EE}{\mathbf{E}}

\def\R{\mathbb{R}}

\def\R{\mathbb R}

\def\EE{\mathbb E}

\def\e{\varepsilon}
\def\la{\langle}
\def\ra{\rangle}
\def\vp{\varphi}
\def\y{\mathbf{y}}

\def\x{\mathbf{x}}

\def\one{{\mathbf 1}}

\def\tf{\tilde{f}}

\def\bld{\boldsymbol}

\newtheorem{lemma}{Lemma}
\newtheorem{theorem}{Theorem}
\newtheorem{definition}{Definition}

\newtheorem{corollary}{Corollary}

\begin{document}
\begin{frontmatter}

\title{Derivative-Free Method For Composite Optimization With Applications To Decentralized Distributed Optimization}


\author[First]{Aleksandr Beznosikov} 
\author[Second]{Eduard Gorbunov} 
\author[Third]{Alexander Gasnikov}

\address[First]{Moscow Institute of Physics and Technology, Russia\\ Sirius University of Science and Technology, Russia\\ (e-mail: beznosikov.an@phystech.edu).}
\address[Second]{Moscow Institute of Physics Technology, Russia\\ Sirius University of Science and Technology, Russia\\ Institute for Information Transmission Problems RAS, Russia \\(e-mail: eduard.gorbunov@phystech.edu).}
\address[Third]{Moscow Institute of Physics and Technology, Russia\\ Sirius University of Science and Technology, Russia \\ Institute for Information Transmission Problems RAS, Russia, \\ Caucasus Mathematical Center, Adyghe State University, Russia \\(e-mail: gasnikov@yandex.ru).}

\begin{abstract}In this paper, we propose a new method based on the Sliding Algorithm from \cite{lan2016gradient,lan} for the convex composite optimization problem that includes two terms: smooth one and non-smooth one. Our method uses the stochastic noised zeroth-order oracle for the non-smooth part and the first-order oracle for the smooth part. To the best of our knowledge, this is the first method in the literature that uses such a mixed oracle for the composite optimization. We prove the convergence rate for the new method that matches the corresponding rate for the first-order method up to a factor proportional to the dimension of the space or, in some cases, its squared logarithm. We apply this method for the decentralized distributed optimization and derive upper bounds for the number of communication rounds for this method that matches known lower bounds. Moreover, our bound for the number of zeroth-order oracle calls per node matches the similar state-of-the-art bound for the first-order decentralized distributed optimization up to to the factor proportional to the dimension of the space or, in some cases, even its squared logarithm.
\end{abstract}


\end{frontmatter}

\section{Introduction}\label{sec:intro}

In this paper we consider finite-sum minimization problem
\begin{equation}
    \min\limits_{x\in X\subseteq \R^n} f(x) = \frac{1}{m}\sum\limits_{i=1}^m f_i(x),\label{eq:finite_sum}
\end{equation}
where each $f_i$ is convex and differentiable function and $X$ is closed and convex. Such kind of problems are highly widespread in machine learning applications \cite{shalev2014understanding}, statistics \cite{spokoiny2012parametric} and control theory \cite{rao2009survey}. In particular, we are interested in the case when functions $f_i$ are stored on different devices which are connected in a network \cite{lan2017communication,scaman2017optimal,scaman2018optimal,scaman2019optimal,9029798,dvinskikh2019decentralized,gorbunov2019optimal,uribe2017optimal}. This scenario often appears when the goal is to accelerate the training of big machine learning models or when the information that defines $f_i$ is known only to the $i$-th \textit{worker}.

In the centralized or parallel case, the general algorithmic scheme can be described in the following way:
\begin{enumerate}
    \item[1)] each worker in parallel performs computations of either gradients or stochastic gradients of $f_i$;
    \item[2)] then workers send the results (not necessarily gradients that they just computed) to one predefined node called \textit{master} node;
    \item[3)] master node processes received information and broadcast new information to each worker that is needed to get new iterate and then the process repeats.
\end{enumerate}
However, such an approach has several problems, e.g.\ synchronization drawback or high requirements to the master node. There are a lot of works that cope with aforementioned drawbacks (see \cite{stich2018local,karimireddy2019error,alistarh2017qsgd,wen2017terngrad}).

Another possible approach to deal with these drawbacks is to use decentralized architecture \cite{bertsekas1989parallel}. Essentially it means that workers are able to communicate only with their neighbors and communications are simultaneous. We also want to mention that such an approach is more robust, e.g.\ it can be applied to time-varying (wireless) communication networks \cite{rogozin2019projected}.

\subsection{Our contributions}
We develop a new method called Zeroth-Order Sliding 
Algorithm ({\tt zoSA}) for solving convex composite problem containing non-smooth part and $L$-smooth part which uses biased stochastic zeroth-order oracle for the non-smooth term and first-order oracle for the smooth component which is, to the best of our knowledge, the first method that uses zeroth-order and first-order oracles for composite optimization problem in such a way (see the details in Section~\ref{sec:main_res}). We prove the convergence result for the proposed method that matches known results for the number of first-oracle calls. 
Regarding the non-smooth component, we prove that the required number of zeroth-order oracle calls is typically $n$ times or, in some cases, $\log^2 n$ times larger then the corresponding bound obtained for the number of first-order oracle calls required for the non-smooth part which is natural for the derivative-free optimization (see \cite{larson2019derivative}). Moreover, we extend the proposed method to the case when the smooth term is additionally strongly convex.

Next, we show how to apply {\tt zoSA} to the decentralized distributed optimization and get results that match the state-of-the-art results for the first-order non-smooth decentralized distributed optimization in terms of the communication rounds.

\section{Notation and Definitions}\label{sec:notation}
We use $\la x,y \ra \eqdef \sum_{i=1}^nx_i y_i$ to denote standard inner product of $x,y\in\R^n$ where $x_i$ corresponds to the $i$-th component of $x$ in the standard basis in $\R^n$. It induces $\ell_2$-norm in $\R^n$ in the following way $\|x\|_2 \eqdef \sqrt{\la x, x \ra}$. We denote $\ell_p$-norms as $\|x\|_p \eqdef \left(\sum_{i=1}^n|x_i|^p\right)^{\nicefrac{1}{p}}$ for $p\in(1,\infty)$ and for $p = \infty$ we use $\|x\|_\infty \eqdef \max_{1\le i\le n}|x_i|$. The dual norm $\|\cdot\|_*$ for the norm $\|\cdot\|$ is defined in the following way: $\|y\|_* \eqdef \max\left\{\la x, y \ra\mid \|x\| \le 1\right\}$.  To denote maximal and minimal positive eigenvalues of positive semidefinite matrix $A\in\R^{n\times n}$ we use $\lambda_{\max}(A)$ and $\lambda_{\min}^{+}(A)$ respectively and we use $\chi(A) \eqdef \nicefrac{\lambda_{\max}(A)}{\lambda_{\min}^{+}(A)}$ to denote condition number of $A$. Operator $\EE[\cdot]$ denotes full mathematical expectation and operator $\EE_\xi[\cdot]$ express conditional mathematical expectation w.r.t. all randomness coming from random variable $\xi$. To define the Kronecker product of two matrices $A\in\R^{m\times m}$ and $B\in\R^{n\times n}$ we use $A \otimes B \in \R^{nm \times nm}$. The identity matrix of the size $n\times n$ is denoted in our paper by $I_n$. 

Since all norms in finite dimensional space are equivalent, there exist such constants $C_1$, $C_2$ and $C_3$ that for all $x\in\R^n$
\begin{eqnarray}
    \|x\|_* \le C_1\|x\|_2,\quad \|x\|_2 \le C_2 \|x\|_*, \quad \|x\| \le C_3\|x\|_2.\label{eq:norm_equiv}
\end{eqnarray}
For example, if $\|\cdot\| = \|\cdot\|_2$, then $C_1 = C_2 = C_3 = 1$ and if $\|\cdot\| = \|\cdot\|_1$, then $\|\cdot\|_* = \|\cdot\|_\infty$ and $C_1 = 1$, $C_2 = C_3 = \sqrt{n}$.

\begin{definition}[$L$-smoothness]
    Function $g$ is called $L$-smooth in $X\subseteq \R^n$ with $L > 0$ w.r.t. norm $\|\cdot\|$ when it is differentiable and its gradient is $L$-Lipschitz continuous in $X$, i.e.\ 
    \begin{equation*}
        \|\nabla g(x) - \nabla g(y)\|_* \le L\|x - y\|,\quad \forall x,y\in X.
    \end{equation*}
\end{definition}
One can show that $L$-smoothness implies (see \cite{nesterov2004introduction})
\begin{equation}
    \label{g-L-smooth} 
    g(x) \leq g(y) + \langle \nabla g(y), x - y\rangle + \dfrac{L}{2}\|x-y\|^2,\quad \forall x,y \in X.
\end{equation}

\begin{definition}[$s$-neighborhood of a set]
    For a given set $X\subseteq \R^n$ and $s > 0$ the $s$-neighborhood of $X$ w.r.t. norm $\|\cdot\|$ is denoted by $X_s$ which is defined as $X_s \eqdef \left\{z\in\R^n\mid \exists\, x\in X:\; \|y-x\| \le s\right\}$. 
\end{definition}

\begin{definition}[Bregman divergence]
    Assume that function $\nu(x)$ is $1$-strongly convex w.r.t. $\|\cdot\|$-norm and differentiable on $X$ function. Then for any two points $x,y\in X$ we define Bregman divergence $V(x,y)$ associated with $\nu(x)$ as follows:
    \begin{equation*}
        V(x,y) \eqdef \nu(y) - \nu(x) - \la\nabla \nu(x), y-x \ra.
    \end{equation*}
\end{definition}
Note that $1$-strong convexity of $\nu(x)$ implies
\begin{equation}
    V(x,y) \ge \frac{1}{2}\|x-y\|^2.\label{eq:bregman_key_property}
\end{equation}
Finally, we denote the Bregman-diameter of the set $X$ w.r.t.\ $V(x,y)$ as $D_{X,V} \eqdef \max\{\sqrt{2V(x,y)}\mid x,y \in X\}$. In view of \eqref{eq:bregman_key_property} $D_{X,V}$ is an upper bound for the standard diameter of the set $D_X \eqdef \max\{\|x-y\|\mid x,y \in X\}$. When $V(x,y) = \frac{1}{2}\|x-y\|_2^2$ (standard Euclidean proximal setup) we have $D_{X,V} = D_X$. If $\|\cdot\| = \|\cdot\|_1$ is $\ell_1$-norm, then in the case when $X$ is a probability simplex, i.e.\ $X = \{x\in \R_+^n\mid \sum_{i=1}^n x_i = 1\}$, and the distance generating function $\nu(x)$ is entropic, i.e.\ $\nu(x) = \sum_{i=1}^n x_i\ln x_i$, we have that $V(x,y)$ is the Kullback-Leibler divergence, i.e.\ $V(x,y) = \sum_{i=1}^n x_i\ln\frac{x_i}{y_i}$, and $D_{X,V} = \sqrt{2\ln{n}}$ (see \cite{ben-tal2015lectures}).

\section{Main Result}\label{sec:main_res}
\subsection{Convex Case}\label{sec:main_res_cvx}
We consider the composite optimization problem
\begin{equation}
    \label{problem_orig} 
    \min\limits_{x \in X} \Psi_0(x) = f(x) + g(x),
\end{equation}
where $X\subseteq \R^n$ is a compact and convex set with diameter $D_X$ in $\|\cdot\|$-norm, function $g$ is convex and $L$-smooth on $X$, $f$ is convex differentiable function on $X$. Assume that we have an access to the first-order oracle for $g$, i.e.\ gradient $\nabla g(x)$ is available, and to the biased stochastic zeroth-order oracle for $f$ (see also \cite{gorbunov2018accelerated}) that for a given point $x$ returns noisy value $\tf(x)$ such that
\begin{equation}
    \label{tilde_f} 
    \tilde{f}(x) \eqdef f(x, \xi) + \Delta(x)
\end{equation}
where $\Delta(x)$ is a bounded noise of unknown nature
\begin{equation}
    \label{delta} 
    |\Delta(x)| \leq \Delta
\end{equation}
and random variable $\xi$ is such that
\begin{equation}
    \label{xi0} 
    \mathbb{E}[f(x,\xi)] = f(x), \\
\end{equation}
Additionally, we assume that for all $x \in X_s$ ($s \le D_X$)
\begin{equation}
    \|\nabla f(x,\xi) \|_2 \leq M(\xi),\quad
    \label{xi} 
    \mathbb{E}[M^2(\xi)] = M^2.
\end{equation}
This assumption implies that for all $x \in X_s$ 
\begin{equation*}
    |f(x,\xi) - f(y, \xi)| \leq M(\xi) \|x-y\|_2
\end{equation*}
and
\begin{equation*}
    \|\nabla f(x)\|_2 \leq M.
\end{equation*}

Using this one can construct a stochastic approximation of $\nabla f(x)$ via finite differences (see \cite{Nesterov,Shamir15}):
\begin{equation}
    \label{oracle_f} 
    \tilde{f}_r'(x) = \frac{n}{2 r} (\tilde{f}(x + r e) - \tilde{f}(x - r e) ) e
\end{equation}
where $e$ is a random vector uniformly distributed on the Euclidean sphere and
\begin{equation}
    r < sC_3 \label{eq:constraint_for_smoothing_parameter}
\end{equation}
is a smoothing parameter. Inequality \eqref{eq:constraint_for_smoothing_parameter} guarantees that the considered approximation requires points only from $s$-neighborhood of $X$ since $\|re\| \le rC_3$ (see \eqref{eq:norm_equiv}). Therefore, throughout the paper we assume that \eqref{eq:constraint_for_smoothing_parameter} holds. Following \cite{Shamir15} we assume that there exists such constant $p_* > 0$ that
\begin{eqnarray}
    \sqrt[4]{\EE[\|e\|_*^4]} &\le& p_*.\label{eq:condition_on_u}
\end{eqnarray}
For example, when $\|\cdot\| = \|\cdot\|_2$ we have $p_* = 1$ and for the case when $\|\cdot\| = \|\cdot\|_1$ one can show that $p_* = O\left(\sqrt{\nicefrac{\ln(n)}{n}}\right)$ (see Corollaries~2~and~3 from \cite{Shamir15}).
Consider also the smoothed version
\begin{equation}
    \label{F} 
    F(x) \eqdef \EE_e [f(x + r e)]
\end{equation}
of $f(x)$ which is a differentiable in $x$ function. In the following we summarize key properties of $F(x)$.
\begin{lemma}[see also Lemma~8 from \cite{Shamir15}]\label{lem:lemma_8_shamir_main}
    Assume that differentiable function $f$ defined on $X_s$ satisfy $\|\nabla f(x)\|_2\le M$ with some constant $M > 0$. Then $F(x)$ defined in \eqref{F} is convex, differentiable and $F(x)$ satisfies
\begin{eqnarray}
    \sup_{x \in X} |F(x) - f(x)| &\leq& rM\label{orig_sm_main},\\    
    \nabla F(x) &=& \mathbb{E}_{e} \left[\frac{n}{r} f(x + r e)e\right],\label{grad_sm_main}\\
    \|\nabla F(x)\|_* &\le& \tilde c p_* \sqrt{n}M,\label{eq:norm_nabla_F_bound_main}
\end{eqnarray}
where $\tilde c$ is some positive constant independent of $n$ and $p_*$ is defined in \eqref{eq:condition_on_u}.
\end{lemma}

In other words, $F(x)$ provides a good approximation of $f(x)$ for small enough $r$. Therefore, instead of solving \eqref{problem_orig} directly one can focus on the problem
\begin{equation}
    \label{problem} 
    \min\limits_{x \in X} \Psi(x) \eqdef F(x) + g(x)
\end{equation}
with small enough $r$ since the difference between optimal values for \eqref{problem_orig} and \eqref{problem} is at most $rM$.
The following lemma establishes useful relations between $\nabla F(x)$ and $\tilde{f}_r'(x)$ defined in \eqref{oracle_f}.
\begin{lemma}[modification of Lemma~10 from \cite{Shamir15}]\label{lem:second_lemma}
For $\tf'_r(x)$ defined in \eqref{oracle_f} the following inequalities hold: 
\begin{equation}
    \|\EE[ \tilde{f}_r'(x)] - \nabla F(x)\|_* \leq  \frac{n\Delta p_*}{r},\label{exp_noise_main}
\end{equation}
\begin{equation}
    \EE[\|\tilde{f}_r'(x)\|^2_{*}] \leq 2p_*^2\left(cnM^2 + \frac{n^2\Delta^2}{r^2}\right)\label{exp_square_main},
\end{equation}
where $c$ is some positive constant independent of $n$.
\end{lemma}
In other words, one can consider $\tilde{f}_r'(x)$ as a biased stochastic gradient of $F(x)$ with bounded second moment and apply Stochastic Gradient Sliding from \cite{lan2016gradient,lan} with this stochastic gradient to solve problem \eqref{problem}.
\begin{algorithm} [H]
	\caption{Zeroth-Order Sliding Algorithm ({\tt zoSA})}
	\label{alg}
	\begin{algorithmic}
\State
\noindent {\bf Input:} Initial point $x_0 \in X$ and iteration limit $N$.
\State Let $\beta_k \in \RR_{++}, \gamma_k \in \RR_+$, and $T_k \in {\mathbb N}$, $k = 1, 2, \ldots$,
be given and
set $\overline x_0 = x_0$. 
\For {$k=1, 2, \ldots, N$ }
    \State 1. Set $\underline x_k = (1 - \gamma_k) \overline x_{k-1} + \gamma_k x_{k-1}$,
    and let $h_k(\cdot) \equiv l_g(\underline x_{k}, \cdot)$ be defined in \eqref{lg}.
    \State 2. Set
    \begin{equation*}
        (x_k, \tilde x_k) = \text{\tt PS}(h_k, x_{k-1}, \beta_k, T_k);
    \end{equation*}
    \State 3. Set $\overline x_k = (1-\gamma_k) \overline x_{k-1} + \gamma_k \tilde x_k$. 
\EndFor
\State 
\noindent {\bf Output:} $\overline x_N$.

\Statex
\Statex The  $\text{\tt{PS} }$(prox-sliding) procedure.

\State {\bf procedure:} {$(x^+, \tilde x^+) = \text{\tt{PS}}$($h$, $x$, $\beta$, $T$)}
\State Let the parameters $p_t \in \R_{++}$ and $\theta_t \in [0,1]$,
$t = 1, \ldots$, be given. Set $u_{0} = \tilde u_0 = x$.
\State {\bf for} $t = 1, 2, \ldots, T$ {\bf do}
    \begin{eqnarray}
        u_{t} &=& \argmin_{u \in X} \Big\{h(u) + \langle \tilde{f}_r'(u_{t-1}), u \rangle\notag\\
        &&\quad+\beta V(x, u) + \beta p_t V(u_{t-1}, u)\Big\},\label{u_t}\\
        \tilde u_t &=&  (1-\theta_t) \tilde u_{t-1} + \theta_t u_t.
        \label{tilde_u_t}
    \end{eqnarray}
\State {\bf end for}
\State Set $x^+ = u_T$ and  $\tilde x^+ = \tilde u_T$.
\State {\bf end procedure:}
\end{algorithmic}
\end{algorithm}
In the Algorithm~\ref{alg} we use the following function
\begin{equation}
    \label{lg} 
    l_g(x,y) \eqdef g(x) + \langle \nabla g(x), y-x \rangle.
\end{equation}
At each iteration of {\tt PS} subroutine the new direction $e$ is sampled independently from previous iterations. We emphasize that we do not need to compute values of $F(x)$ which in the general case requires numerical computation of integrals over a sphere. In contrast, our method requires to know only noisy values of $f$ defined in \eqref{tilde_f}.

Next, we present the convergence analysis of {\tt zoSA} that relies on the analysis for the Gradient Sliding method from \cite{lan2016gradient,lan}. The following lemma provides an analysis of the subroutine {\tt PS} from Algorithm~\ref{alg}.

\begin{lemma}[modification of Proposition~8.3 from \cite{lan}]\label{lem:third_lemma}
    Assume that $\{p_t\}_{t\ge 1}$ and $\{\theta_t\}_{t\ge 1}$ in the subroutine {\tt PS} of Algorithm~\ref{alg} satisfy
    \begin{eqnarray}
        \theta_t &=& \frac{P_{t-1} - P_t}{(1 - P_t)P_{t-1}},\label{p_t_main}\\
        P_t &=& 
        \begin{cases}
            1 & t = 0,\\
            p_t(1 + p_t)^{-1}P_{t-1} & t \geq 1.
        \end{cases}\notag
    \end{eqnarray}
    Then for any $t \geq 1$ and $u\in X$:
    \begin{eqnarray}
        \label{lemma_2_main}
        \beta (1 - P_t)^{-1}V(u_t, u) + \left[\Phi(\tilde u_t) - \Phi(u)\right] &\notag\\
        &\hspace{-7.5cm}\leq \beta P_t(1 - P_t)^{-1}V(u_0,u)\notag\\
        &\hspace{-5cm}+ P_t(1 - P_t)^{-1} \sum\limits_{i=1}^t (p_i P_{i-1})^{-1} \Big[ \frac{(\tilde M + \| \delta_i\|_*)^2}{2 \beta p_i}\notag\\
        &\hspace{-1.5cm}+ \langle \delta_i, u-u_{i-1} \rangle\Big],
    \end{eqnarray}
    where 
    \begin{equation}
        \label{Phi_main}
        \Phi(u) = h(u) + F(u) + \beta V(x,u),
    \end{equation}
    \begin{equation}
        \label{delta_main}
        \delta_t = \tilde f_r' (u_{t-1}) - \nabla F(u_{t-1}).
    \end{equation}
    \begin{eqnarray*}
        \tilde M = c \sqrt{n}C_1M,
    \end{eqnarray*}
    $c\;$ is some positive constant independent of $n$, $C_1$ is from \eqref{eq:norm_equiv}.
\end{lemma}
Using the lemma above we derive the main result.
\begin{theorem}\label{thm:first_theorem} Assume that $\{ p_t\}_{t\ge 1}$, $\{\theta_t\}_{t\ge 1}$, $\{\beta_k\}_{k \ge 1}$, $\{\gamma_k\}_{k\ge 1}$ in Algorithm~\ref{alg} satisfy \eqref{p_t_main} and
\begin{equation}
    \label{gamma_main}
    \gamma_1 = 1,~~~~\beta_k - L \gamma_k \geq 0,~~~ k\geq 1,
\end{equation}
\begin{equation}
    \label{gamma_k_main}
    \frac{\gamma_k \beta_k}{\Gamma_k(1 - P_{T_k})} \leq \frac{\gamma_{k-1} \beta_{k-1}}{\Gamma_{k-1}(1 - P_{T_{k-1}})} ,~~~ k\geq 2.
\end{equation}
Then 
\begin{eqnarray}
    \label{t2_1}
    \EE[\Psi(\overline x_N) - \Psi(x^*)] &\notag\\
    &\hspace{-3cm}\leq \frac{\Gamma_N \beta_1}{1 - P_{T_1}} V(x_0,u) +
    \Gamma_N \sum\limits_{k=1}^N\sum\limits_{i=1}^{T_k} \Bigg[ \frac{(\tilde M^2 + \sigma^2)\gamma_k P_{T_k}}{\beta_k \Gamma_k(1 - P_{T_k})p_i^2 P_{i-1}}
    \notag\\ 
    &\hspace{-0.5cm}+ \frac{n \Delta D_X p_*}{r}\cdot\frac{\gamma_k P_{T_k}}{\Gamma_k(1 - P_{T_k})p_i P_{i-1}} \Bigg],
\end{eqnarray}
where $x^*$ is an arbitrary optimal point for \eqref{problem}, $P_t$ is from \eqref{p_t_main}, 
\begin{equation} \Gamma_k = 
\begin{cases} 
            1, & k = 1,\\
            (1 - \gamma_k)\Gamma_{k-1}, & k > 1
        \label{gamma_kk_main}
\end{cases}
\end{equation}
and 
\begin{equation}
\sigma^2 \eqdef 4p_*^2\left(CnM^2 + \frac{n^2\Delta^2}{r^2} \right),
\label{sig_main}
\end{equation}
where $C$ is some positive constant independent of $n$.
\end{theorem}
    
The next corollary suggests the particular choice of parameters and states convergence guarantees in a more explicit way.
\begin{corollary}\label{cor:main} Suppose that $\{ p_t\}_{t\ge 1}$, $\{\theta_t\}_{t\ge 1}$ are
\begin{equation}
    \label{pt_main}
    p_t = \frac{t}{2}, ~~~ \theta_t = \frac{2(t+1)}{t(t+3)}, ~~~\forall t \geq 1,
\end{equation}
$N$ is given, $\{\beta_k \}$, $\{\gamma_k\}$, $T_k$ are
\begin{equation}
    \label{k_main}
    \beta_k = \frac{2L}{k},~~~\gamma_k = \frac{2}{k+1},~~~T_k = \frac{N(\tilde M^2 + \sigma^2)k^2}{\tilde D L^2}
\end{equation}
for $\tilde D = \nicefrac{3D_{X,V}^2}{4}$. Then $\forall N \geq 1$
\begin{equation}
    \label{col1_main}
    \mathbb{E}[\Psi(\overline x_N)- \Psi(x^*)] \leq\frac{12LD_{X,V}^2}{N(N+1)}+ \frac{n \Delta D_Xp_*}{r}.
\end{equation}
\end{corollary}
Finally, we extend the result above to the initial problem \eqref{problem_orig}.
\begin{corollary}\label{cor:main_corollary} Under the assumptions of Corollary~\ref{cor:main} we have that the following inequality holds for all $N \ge 1$:
    \begin{eqnarray}
        \mathbb{E}[\Psi_0(\overline x_N)- \Psi_0(x^*)] &\leq& 2rM + \frac{12LD_{X,V}^2}{N(N+1)}\notag\\
        &&\quad+ \frac{n \Delta D_Xp_*}{r}.\label{original_main}
    \end{eqnarray}
From \eqref{original_main} it follows that if
\begin{eqnarray}
   r &=& \Theta\left(\frac{\varepsilon}{M}\right),\quad \Delta = O\left(\frac{\varepsilon^2}{nMD_X\min\{p_*,1\}}\right)\label{r_delta_main}
\end{eqnarray}
and $\varepsilon = O\left(\sqrt{n}MD_X\right)$, $s = \Omega\left(\nicefrac{\e}{MC_3}\right)$, then the number of evaluations for $\nabla g$ and $\tilde f'_r$, respectively, required by Algorithm~\ref{alg} to find an $\varepsilon$-solution of \eqref{problem_orig}, i.e.\ such $\overline x_N$ that $\EE[\Psi_0(\overline x_N)] - \Psi_0(x^*) \le \varepsilon$, can be bounded by
\begin{equation}
    \label{bound_out_main}
    O\left(\sqrt{\frac{L D_{X,V}^2}{\varepsilon}} \right),
    \end{equation}
\begin{equation}
    \label{bound_in_main}
  O\left(\sqrt{\frac{L D_{X,V}^2}{\varepsilon}} + \frac{D_{X,V}^2nM^2(C_1^2 +p_*^2)}{\varepsilon^2}\right).
    \end{equation}
\end{corollary}


Let us discuss the obtained result and especially bounds \eqref{bound_out_main} and \eqref{bound_in_main}. First of all, consider Euclidean proximal setup, i.e. $\|\cdot\| = \|\cdot\|_2$, $V(x,y) = \frac{1}{2}\|x-y\|_2^2,$ $D_{X,V} = D_X$. In this case we have $p_* = C_1 = C_2 = C_3 = 1$ and bound \eqref{bound_in_main} for the number of \eqref{tilde_f} oracle calls reduces to
\begin{equation*}
  O\left(\sqrt{\frac{L D_X^2}{\varepsilon}} + \frac{D_X^2nM^2}{\varepsilon^2}\right)
\end{equation*}
and the number of $\nabla g(x)$ computations remains the same. It means that our result gives the same number of first-order oracle calls as in the original Gradient Sliding algorithm, while the number of the biased stochastic zeroth-order oracle calls is $n$ times larger in the leading term than in the analogous bound from the original first-order method. In the Euclidean case our bounds reflect the classical dimension dependence for the derivative-free optimization (see \cite{larson2019derivative}).

Secondly, we consider the case when $X$ is the probability simplex in $\R^n$ and the proximal setup is entropic (see the end of Section~\ref{sec:notation}). As we mentioned earlier in Section~\ref{sec:notation} and in the beginning of this section, in this situation we have $D_{X,V} = \sqrt{2\ln n}$, $D_X = 2$, $p_* = O\left(\nicefrac{\ln(n)}{n}\right)$ and $C_1 = 1$, $C_2 = C_3 = \sqrt{n}$. Then number of $\nabla g(x)$ calculations is bounded by $O\left( \sqrt{\nicefrac{(L\ln^2 n)}{\e}}\right)$. As for the number of $\tilde{f}'_r(x)$ computations, we get the following bound:
\begin{equation*}
  O\left(\sqrt{\frac{L \ln^2 n}{\varepsilon}} + \frac{M^2\ln^2 n}{\varepsilon^2}\right).
\end{equation*}
Clearly, in this case we have only polylogarithmical dependence on the dimension instead.

\subsection{Strongly Convex Case}\label{sec:main_res_str_cvx}
In this section we additionally assume that $g$ is $\mu$-strongly convex w.r.t.\ Bregman divergence $V(x,y)$, i.e.\ $\forall x,y\in X$ 
\begin{equation*}
    g(x) \geq g(y) + \langle \nabla g(y), x - y\rangle + \mu V(x,y).
\end{equation*}

Similarly to the original work \cite{lan2016gradient} we use restarts technique in this case and get Algorithm~\ref{alg2}. 
\begin{algorithm} [H]
	\caption{The Multi-phase Zeroth-Order Sliding Algorithm ({\tt M-zoSA})}
	\label{alg2}
	\begin{algorithmic}
\State
{\bf Input:} Initial point $y_0 \in X$ and iteration limit $N_0$, initial estimate $\rho_0$ (s.t. $\Psi (y_0) - \Psi^* \leq \rho_0$)
\For {$i = 1, 2, \ldots, I$ }
    \State Run {\tt zoSA} with $x_0 = y_{i-1}$, $N = N_0$, $\{p_t\}$ and $\{\theta_t\}$ in \eqref{pt_main}, $\{\beta_k\}$ and $\{\gamma_k\}$, $\{T_k\}$ in \eqref{k_main} with $\tilde D = \nicefrac{\rho_0}{\mu 2^i}$, and $y_i$ is output.
\EndFor
\State {\bf Output:} $y_I$.
	\end{algorithmic}
\end{algorithm}

The following theorem states the main complexity results for {\tt M-zoSA}.
\begin{theorem}\label{thm:mzosa_main}
    For {\tt M-zoSA} with $N_0 = 2 \lceil \sqrt{\nicefrac{5L}{\mu}}\rceil$ we have
    \begin{equation}
        \label{converg_str_conv}
        \EE{[\Psi (y_i) - \Psi(y^*)]} \leq \frac{\rho_0}{2^i} + \frac{2 n \Delta D_X p_*}{r}.
    \end{equation}
\end{theorem}
Using this we derive the complexity bounds for {\tt M-zoSA}.
\begin{corollary}\label{cor:mzosa_main}
For all $N \ge 1$ the iterates of {\tt M-zoSA} satisfy
    \begin{equation}
        \label{original_strconv}
        \mathbb{E}[\Psi_0(y_i)- \Psi_0(y^*)] \leq 2rM + \frac{\rho_0}{2^i} + \frac{2n \Delta D_Xp_*}{r}.
    \end{equation}
From \eqref{original_strconv} it follows that if
\begin{eqnarray}
    r &=& \Theta\left(\frac{\varepsilon}{M}\right),\quad \Delta = O\left(\frac{\varepsilon^2}{nMD_X\min\{p_*,1\}}\right)\label{r_delta_strconv}
\end{eqnarray}
and $\varepsilon = O\left(\sqrt{n}MD_X\right)$, $s = \Omega\left(\nicefrac{\e}{MC_3}\right)$, then the number of evaluations for $\nabla g$ and $\tilde f'_r$, respectively, required by Algorithm~\ref{alg2} to find a $\varepsilon$-solution of \eqref{problem_orig} can be bounded by
\begin{equation}
    \label{bound_out_stconv}
    O\left(\sqrt{\frac{L}{\mu}} \log_2 \max\left[ 1, \nicefrac{\rho_0}{\varepsilon}\right]\right),
    \end{equation}
\begin{equation}
    \label{bound_in_stconv}
   O\left(\sqrt{\frac{L}{\mu}} \log_2 \max\left[ 1, \nicefrac{\rho_0}{\varepsilon}\right] + \frac{nM^2(C_1^2 +p_*^2)}{\mu\varepsilon}\right).
    \end{equation}
\end{corollary}

\section{From Composite Optimization to Convex Optimization with Affine Constraints and Decentralized Distributed Optimization}\label{sec:affine_and_decentral}
In this section we apply the obtained results to the convex optimization problems with affine constraints and after that to the decentralized distributed optimization problem.
\subsection{Convex Optimization with Affine Constraints}\label{sec:affine}
As an intermediate step between the composite optimization problem \eqref{problem_orig} and decentralized distributed optimization we consider the following problem
\begin{equation}
\label{PP}
\min_{Ax=0, x\in X} f(x),    
\end{equation}
where $A \succeq 0$ and $\text{Ker} A \neq \{0\}$ and $X$ is convex compact in $\R^n$ with diameter $D_X$. The dual problem for \eqref{PP} can be written in the following way
\begin{eqnarray}
 && \min_{y} \psi(y),\quad \text{where}\label{DP}\\
\varphi(y) &=& \max_{x\in X}\left\{\langle y,x\rangle - f(x)\right\},\nonumber\\
\psi(y) &=& \vp(A^\top y) = \max_{x\in Q}\left\{\langle y,Ax\rangle - f(x)\right\}\nonumber\\
&=& \langle y,Ax(A^\top y)\rangle - f(x(A^\top y))\notag\\
&=&  \langle A^\top y,x(A^\top y)\rangle - f(x(A^\top y)),\notag
\end{eqnarray}
where $x(y) \eqdef \argmax_{x\in X}\left\{\la y, x\ra - f(x)\right\}$. The solution of \eqref{DP} with the smallest $\ell_2$-norm is denoted in this paper as $y_*$. This norm $R_y \eqdef \|y_*\|_2$ can be bounded as follows \cite{lan2017communication}:
\begin{equation*}
R_y^2 \le \frac{\|\nabla f(x^*)\|_2^2}{\lambda_{\min}^{+}(A^\top A)}.
\end{equation*}

Following \cite{gasnikov2018universalgrad,dvinskikh2019decentralized,gorbunov2019optimal} we consider the penalized problem
\begin{equation}
\label{penalty}
\min_{x\in X} F(x) = f(x) + \frac{R_y^2}{\e}\| Ax\|_2^2,  
\end{equation}
where $\e > 0$ is some positive number. It turns out (see the details in \cite{gorbunov2019optimal}) that if we have such $\hat{x}$ that $F(\hat x) - \min_{x\in X} F(x) \le \e$ then we also have
\begin{equation*}
    f(\hat x) - \min\limits_{Ax = 0, x\in X} f(x) \le \e,\quad \|A\hat x\|_2 \le \frac{2\e}{R_y}.
\end{equation*}
We notice that this result can be generalized in the following way: if we have such $\hat{x}$ that $\EE[F(\hat x)] - \min_{x\in X} F(x) \le \e$ then we also have
\begin{equation}
    \EE[f(\hat x)] - \min\limits_{Ax = 0, x\in X} f(x) \le \e,\quad \sqrt{\EE[\|A\hat x\|_2^2]} \le \frac{2\e}{R_y}.\label{eq:F(x^N)_guarantee_consequence_stoch}
\end{equation}

Next, we consider the problem \eqref{penalty} as \eqref{problem_orig} with $g(x) = \nicefrac{R_y^2\|Ax\|_2^2}{\e}$. Assume that $\|\nabla f(x)\|_2 \le M$ for all $x\in X$ and for $f$ we have an access to the biased stochastic oracle defined in \eqref{tilde_f}. We are interested in the situation when $\nabla g(x) = \nicefrac{2R_y^2 A^\top Ax}{\e}$ can be computed exactly. Moreover, it is easy to see that $g(x)$ is $\nicefrac{2R_y^2 \lambda_{\max}(A^\top A)}{\e}$-smooth w.r.t. $\ell_2$-norm. Applying Corollary~\ref{cor:main_corollary} we get that in order to produce such a point $\hat{x}$ that satisfies \eqref{eq:F(x^N)_guarantee_consequence_stoch} Algorithm~\ref{alg} applied to solve \eqref{penalty} requires
\begin{equation*}
    O\left(\sqrt{\frac{\lambda_{\max}(A^\top A) R_y^2 D_X^2}{\varepsilon^2}} \right) \text{ calculations of $A^\top Ax$} 
\end{equation*}
and
\begin{equation*}
    O\left(\sqrt{\frac{\lambda_{\max}(A^\top A) R_y^2 D_X^2}{\varepsilon^2}} + \frac{nD_X^2M^2}{\varepsilon^2}\right)
\end{equation*}
calculations of $\tf(x)$ since $p_* = C_2 = C_1 = 1$ for the Euclidean case. As we mentioned at the end of Section~\ref{sec:main_res}, this bound depends on dimension $n$ in the classical way.

\subsection{Decentralized Distributed Optimization}\label{sec:decentral}
Now, we go back to the problem \eqref{eq:finite_sum} and, following \cite{scaman2017optimal}, we rewrite it in the distributed fashion:
\begin{equation}
    \min\limits_{\substack{x_1 = \ldots = x_m\\ x_1,\ldots,x_m \in X}} f(\x) = \frac{1}{m}\sum\limits_{i=1}^m f_i(x_i),\label{eq:distrib_rewritten}
\end{equation}
where $\x^\top = (x_1^\top,\ldots, x_m^\top)^\top \in \R^{nm}$. Recall that we consider the situation when $f_i$ is stored on the $i$-th node. In this case one can interpret $x_i$ from \eqref{eq:distrib_rewritten} as a local variable of $i$-th node and $x_1 = \ldots = x_n$ as a consensus condition for the network. The common trick \cite{scaman2017optimal,scaman2018optimal,scaman2019optimal,uribe2017optimal} to handle this condition is to rewrite it using the notion of Laplacian matrix. In general, the Laplacian matrix $\overline{W} = \|\overline{W}_{ij}\|_{i,j=1,1}^{m,m} \in \R^{m\times m}$ of the graph $G$ with vertices $V$, $|V| = m$ and edges $V$ is defined as follows: \begin{equation*}
    \overline{W}_{ij} = \begin{cases} 
    -1, &\text{if } (i,j)\in E,\\
    \deg(i), &\text{if } i=j,\\
    0 &\text{otherwise},
    \end{cases}
\end{equation*}
where $\deg(i)$ is degree of $i$-th node. In this paper we focus only on the connected networks. In this case $\overline{W}$ has unique eigenvector $\bld{1}_m \eqdef (1,\ldots,1)^\top \in \R^m$ associated to the eigenvalue $0$. Using this one can show that for all vectors $a = (a_1,\ldots,a_m)^\top\in \R^m$ we have the following equivalence:
\begin{equation}
    a_1 = \ldots = a_m \; \Longleftrightarrow \; Wa = 0.\label{eq:main_property_of_laplacian_simple}
\end{equation}
Using the Kronecker product $W \eqdef \overline{W}\otimes I_n$, which is also called Laplacian matrix for simplicity, one can generalize \eqref{eq:main_property_of_laplacian_simple} for the $n$-dimensional case:
\begin{equation*}
    x_1 = \ldots = x_m \; \Longleftrightarrow \; W\x = 0
\end{equation*}
and
\begin{equation*}
    x_1 = \ldots = x_m \; \Longleftrightarrow \; \sqrt{W}\x = 0.
\end{equation*}
That is, instead of the problem \eqref{eq:distrib_rewritten} one can consider the equivalent problem
\begin{equation}
    \min\limits_{\substack{\sqrt{W}\x = 0, \\ x_1,\ldots, x_m \in X}} f(\x) = \frac{1}{m}\sum\limits_{i=1}^m f_i(x_i).\label{eq:main_problem_decentralized_sec_rewritten}
\end{equation}

Next, we need to define parameters of $f$ using local parameters of $f_i$. Assume that for each $f_i$ we have $\|f_i(x_i)\|_2 \le M$ for all $x_i \in X$, all $f_i$ are convex functions, the starting point is $\x_0^\top = (x_0^\top,\ldots,x_0^\top)^\top$ and $\x_*^\top = (x_*^\top,\ldots, x_*^\top)^\top$ is the optimality point for \eqref{eq:main_problem_decentralized_sec_rewritten}. Then, one can show (see \cite{gorbunov2019optimal} for the details) that $\|\nabla f(\x)\|_2 \le \nicefrac{M}{\sqrt{m}}$ on the set of such $\x$ that $x_1,\ldots, x_m\in X$, $D_{X^m}^2 = mD_X^2$ and $R_{\y}^2 \eqdef \|\y_*\|_2^2 \le \nicefrac{M^2}{m\lambda_{\min}^+(W)}.$

Now we are prepared to apply results obtained in Section~\ref{sec:affine} to the problem \eqref{eq:main_problem_decentralized_sec_rewritten}. Indeed, this problem can be viewed as \eqref{eq:main_problem_decentralized_sec_rewritten} with $A = \sqrt{W}$. Taking this into account, we conclude that one $A^\top Ax$ calculation corresponds to the calculation of $Wx$ which can be computed during one communication round in the network with Laplacian matrix $W$. This simple observation implies that in order to produce such a point $\hat \x$ that satisfies \eqref{eq:F(x^N)_guarantee_consequence_stoch} with $\hat x = \hat \x$, $A := \sqrt{W}$, $X := X^n$, $R_y := R_{\y}$ Algorithm~\ref{alg} applied to the penalized problem \eqref{penalty} requires 
\begin{equation*}
    O\left(\sqrt{\frac{\chi(W) M^2 D_X^2}{\varepsilon^2}} \right) \text{ communication rounds}
\end{equation*}
and
\begin{equation*}
    O\left(\sqrt{\frac{\chi(W) M^2 D_X^2}{\varepsilon^2}} +  \frac{nD_X^2M^2}{\varepsilon^2}\right) 
\end{equation*}
calculations of $\tf(x)$ per node since $p_* = 1$ for the Euclidean case. The bound for the communication rounds matches the lower bound from \cite{scaman2018optimal,scaman2019optimal} and we conjecture that under our assumptions the obtained bound for zeroth-order oracle calculations per node is optimal up to polylogarithmic factors in the class of methods with optimal number of communication rounds (see also \cite{dvinskikh2019decentralized,gorbunov2019optimal}).

\section{Discussion}


To conclude, the proposed method~--- {\tt zoSA}~--- is the first, to the best of our knowledge, $\nicefrac{1}{2}$-order method for the convex composite optimization: it uses zeroth-order oracle for the non-smooth term and the first-order oracle for the smooth one. It has solid theory and is competitive in practice even with some first-order methods (see our numerical experiments in the appendix).

As for the future work, it would be interesting to study zeroth-order distributed methods for the smooth decentralized distributed optimization using the technique from \cite{gorbunov2019optimal}. Another direction for future research is in developing the analysis of the proposed method for the case when $X$ is unbounded and, in particular, when $X = \R^n$ via recurrences techniques from \cite{gorbunov2018accelerated,gorbunov2019optimal}.

\begin{ack}
The research of A.~Beznosikov, E.~Gorbunov and A.~Gasnikov was partially supported by RFBR, project number 19-31-51001.
The research of E.~Gorbunov was also partially supported by the Ministry of Science and Higher Education of the Russian Federation (Goszadaniye) 075-00337-20-03 and the research of A.~Gasnikov was also partially supported by Yahoo! Research Faculty Engagement Program.

\end{ack}

{ \bibliography{ifacconf}}             
                                                   






\newpage
\onecolumn
\part*{Appendix. Derivative-Free Method For Composite Optimization With Applications To Decentralized Distributed Optimization}

\section{Numerical Experiments}\label{sec:experiments}
\begin{figure*}[h!]
\centering
\includegraphics[width =  0.24\textwidth]{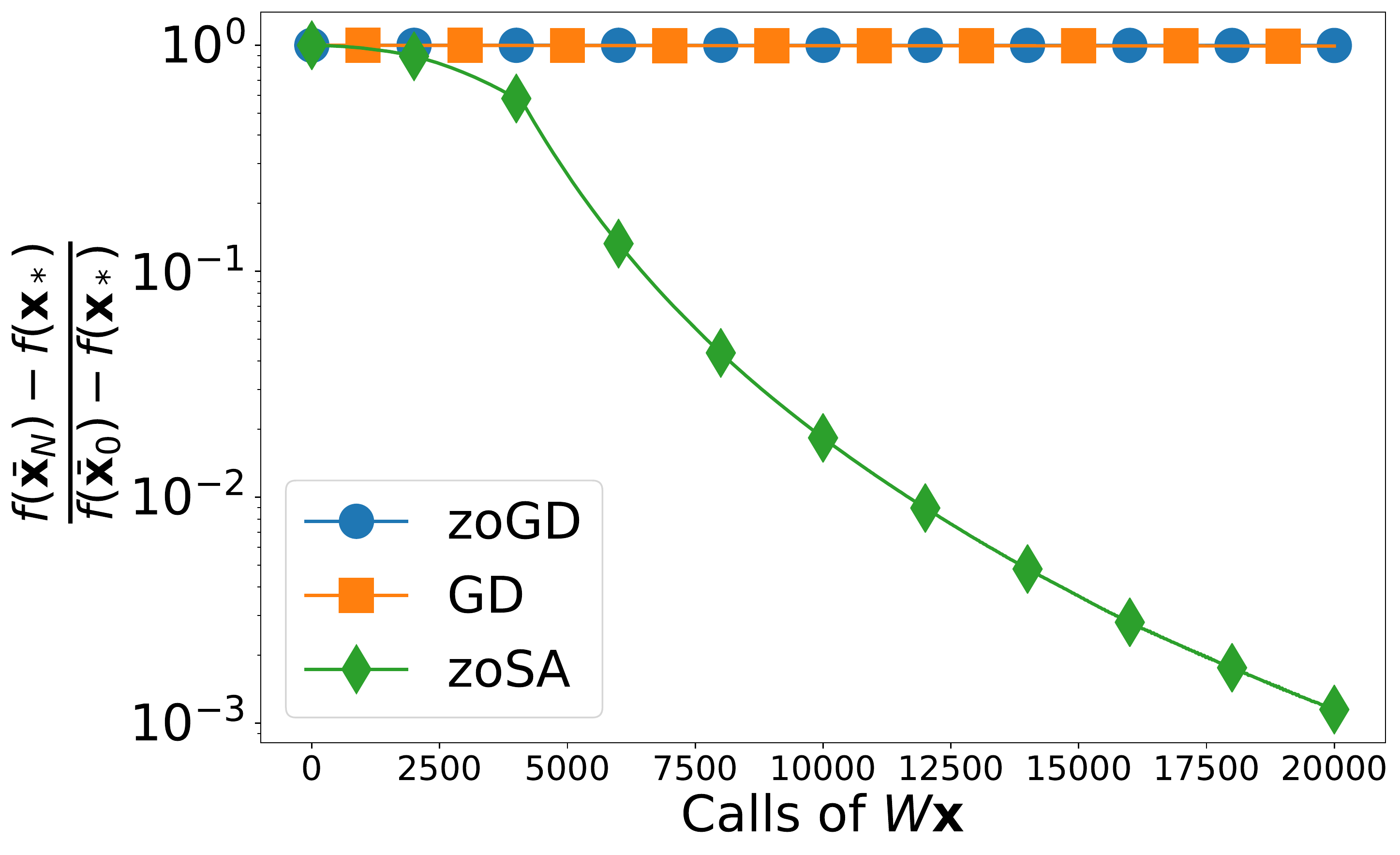}
\includegraphics[width =  0.24\textwidth]{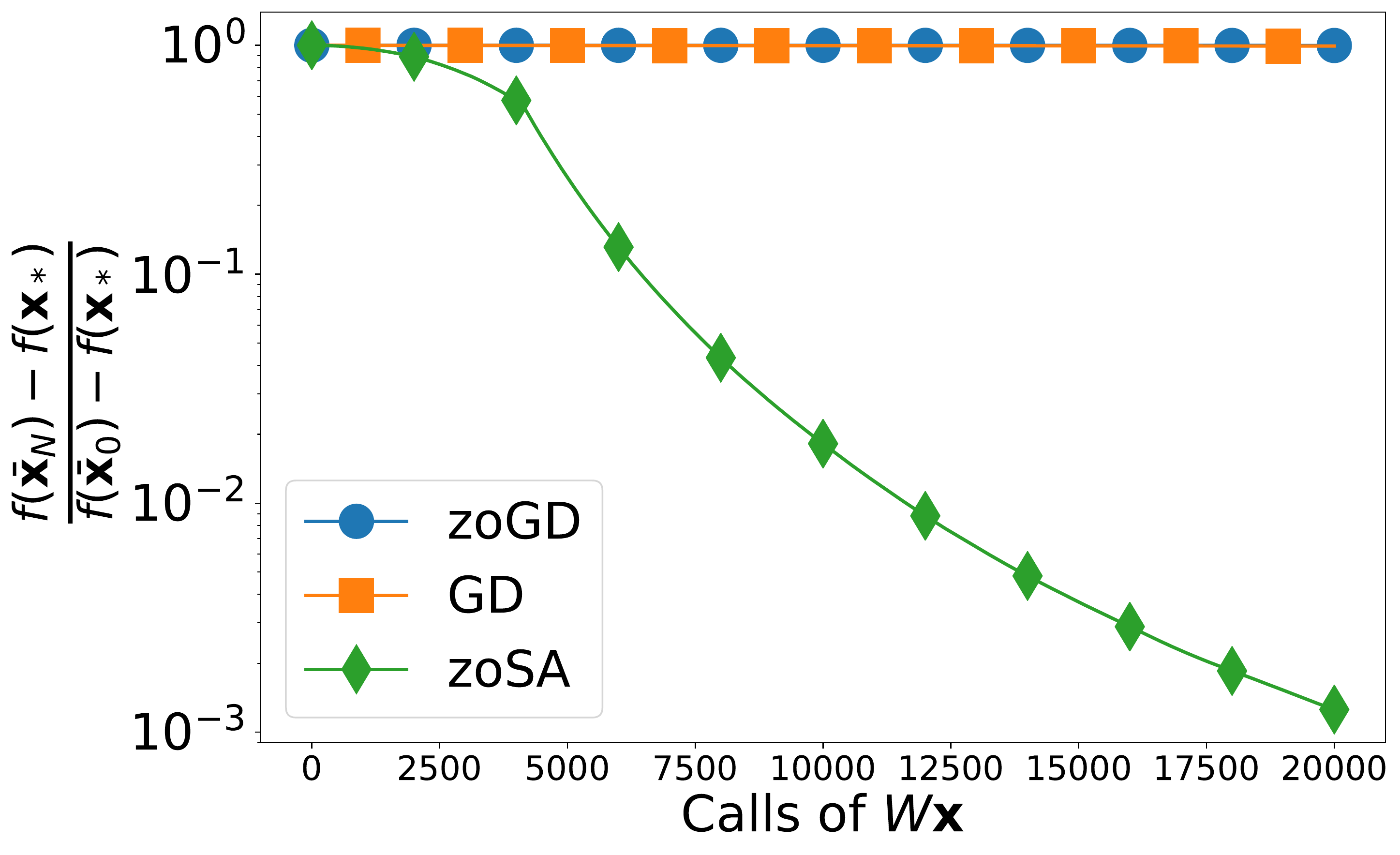}
\includegraphics[width =  0.24\textwidth]{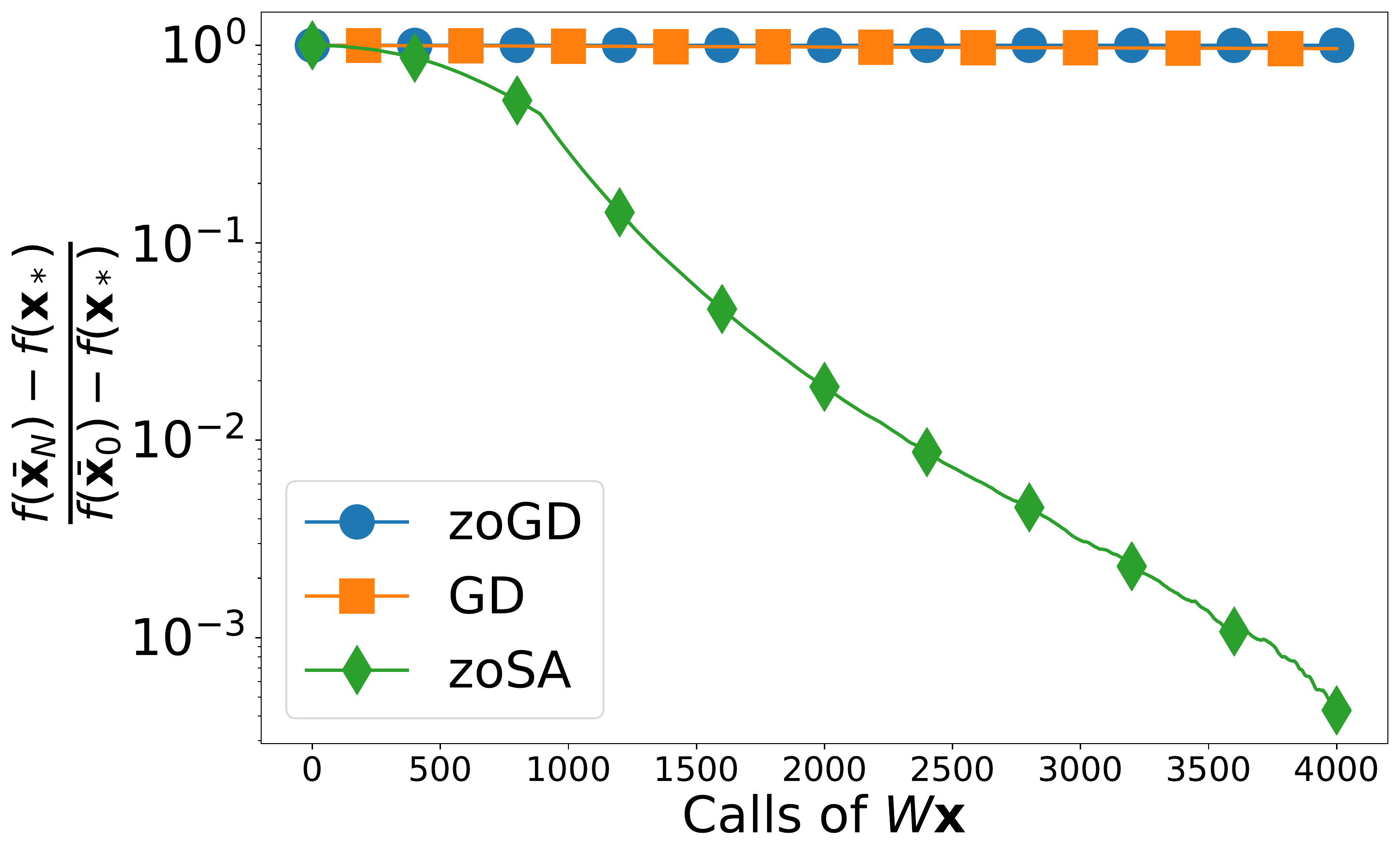}
\includegraphics[width =  0.24\textwidth]{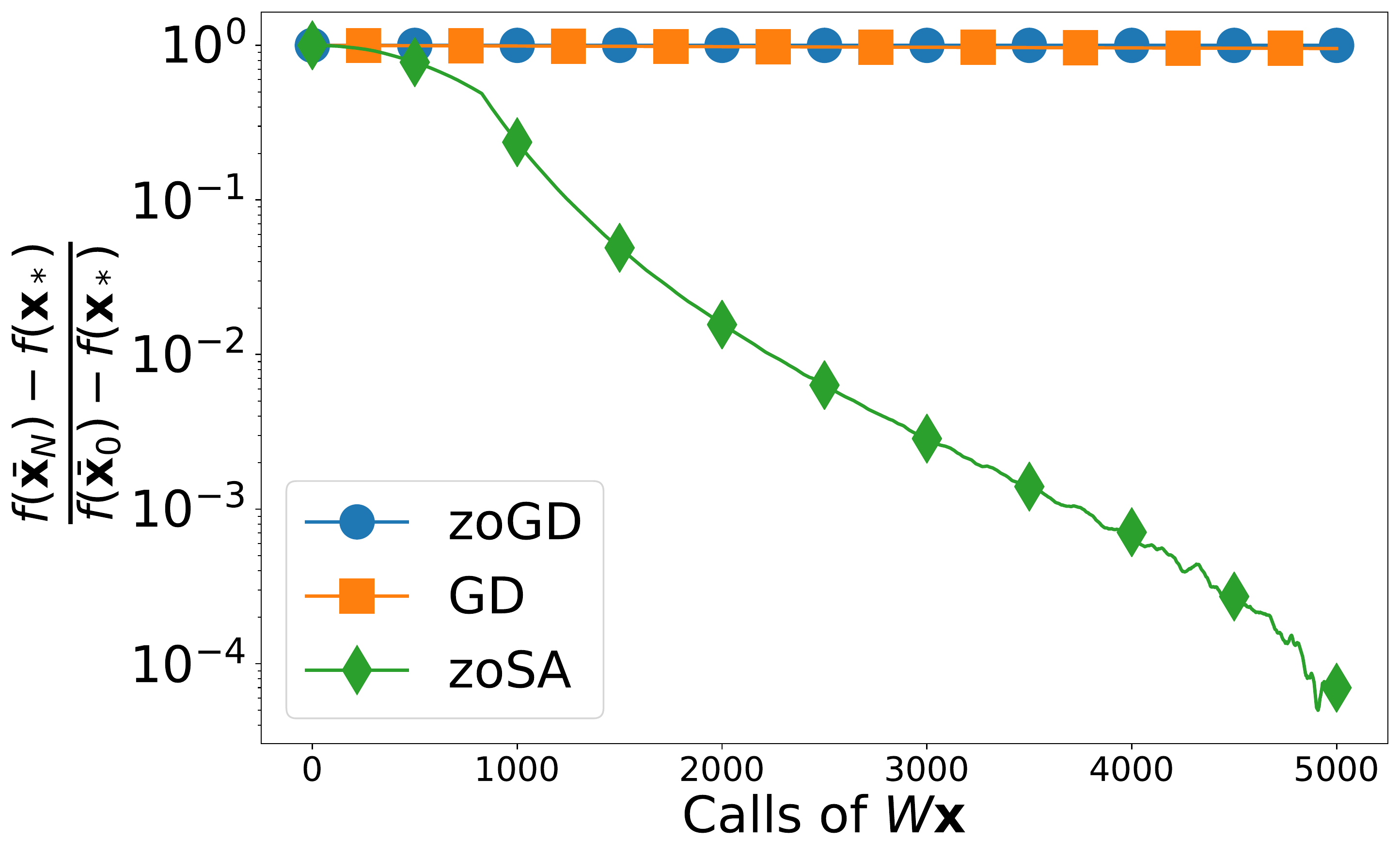}
\includegraphics[width =  0.24\textwidth]{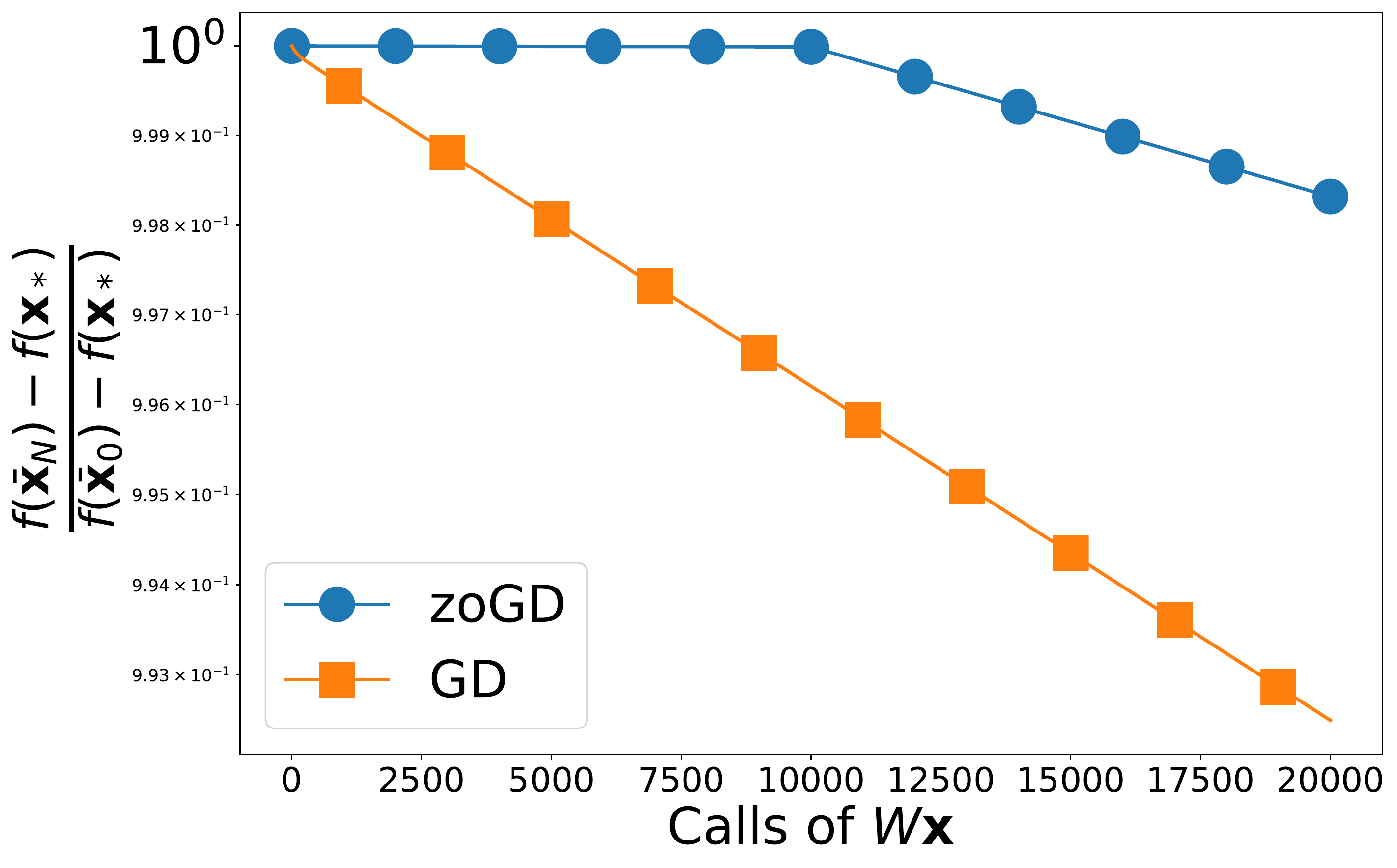}
\includegraphics[width =  0.24\textwidth]{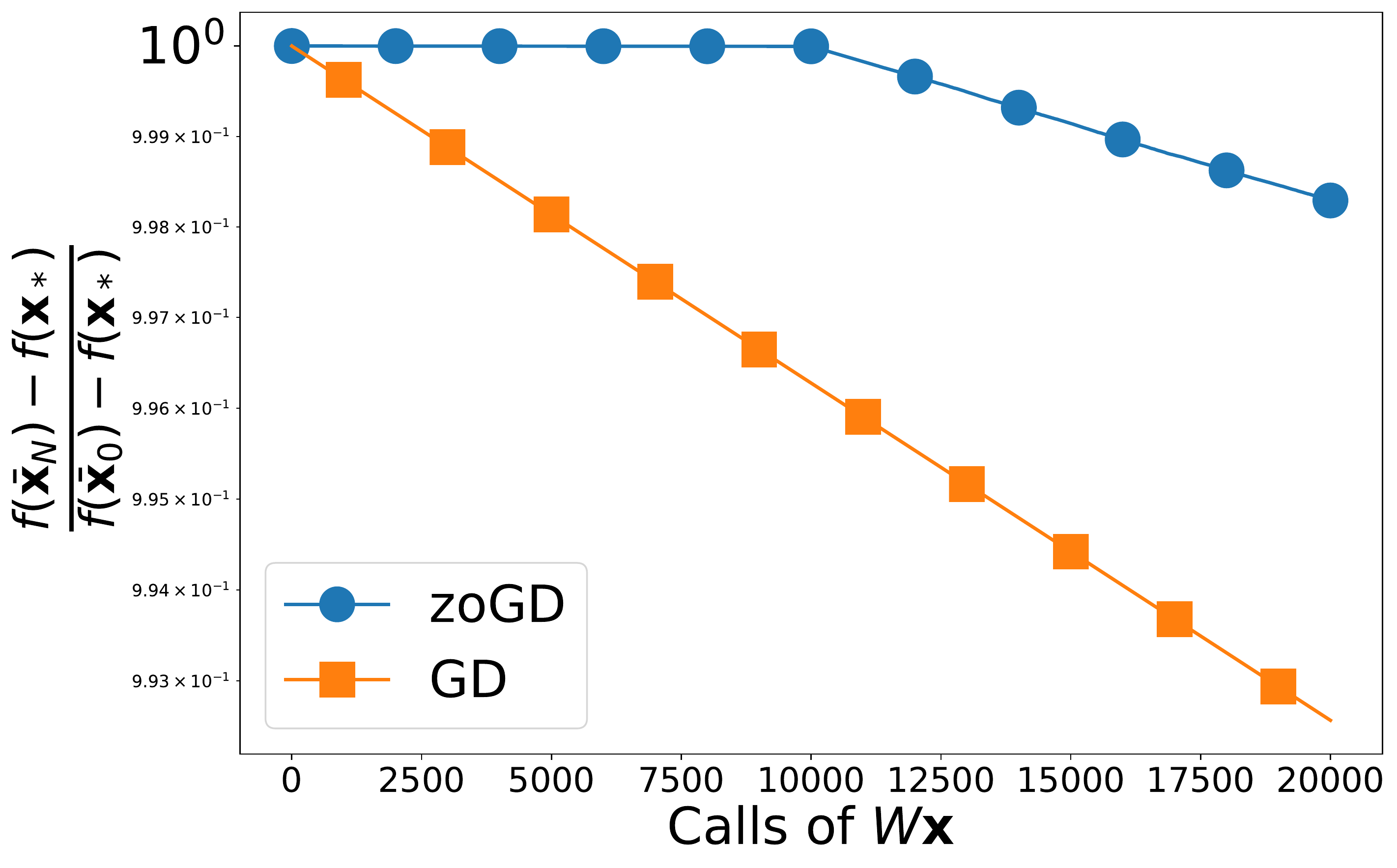}
\includegraphics[width =  0.24\textwidth]{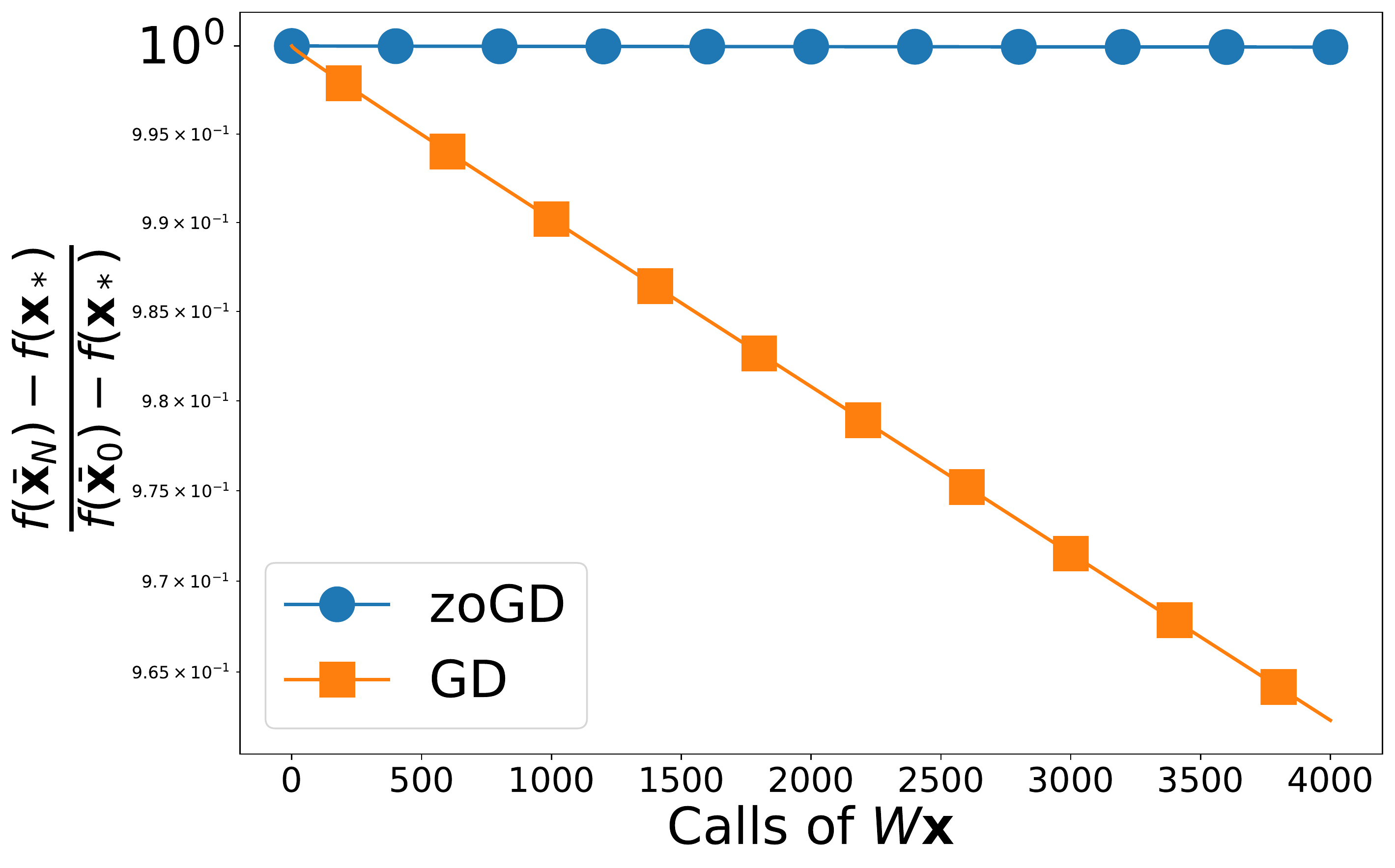}
\includegraphics[width =  0.24\textwidth]{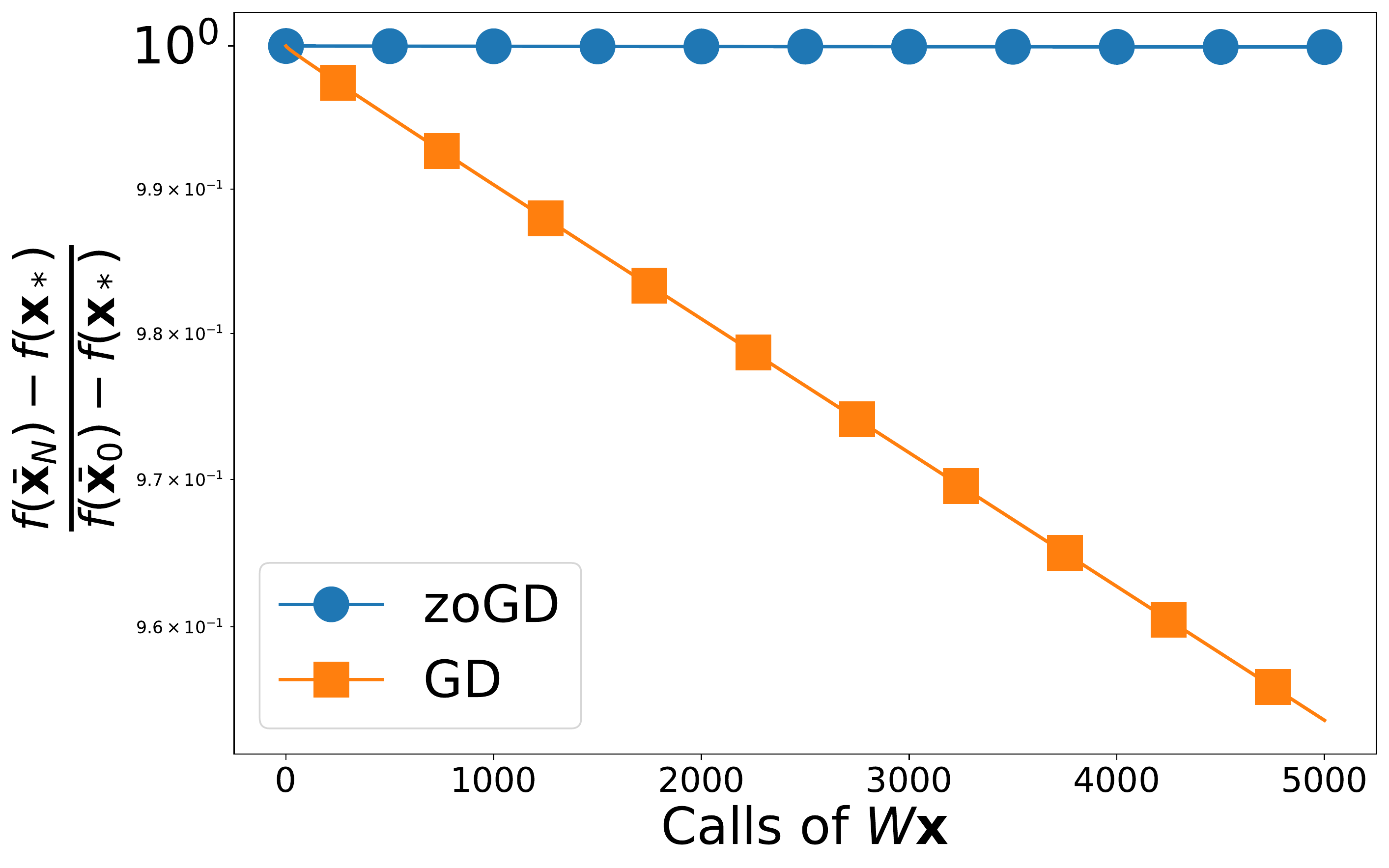}
\includegraphics[width =  0.24\textwidth]{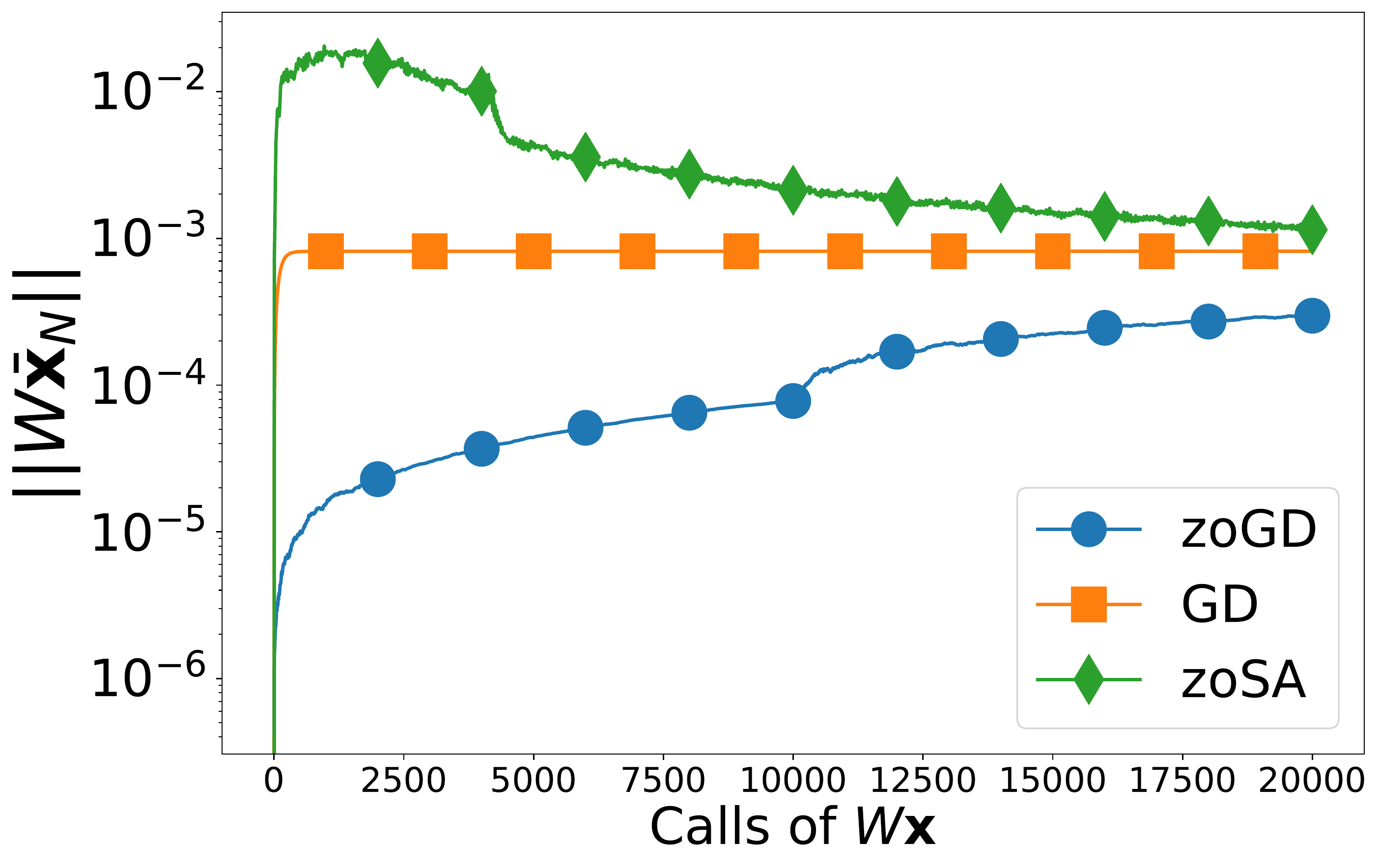}
\includegraphics[width =  0.24\textwidth]{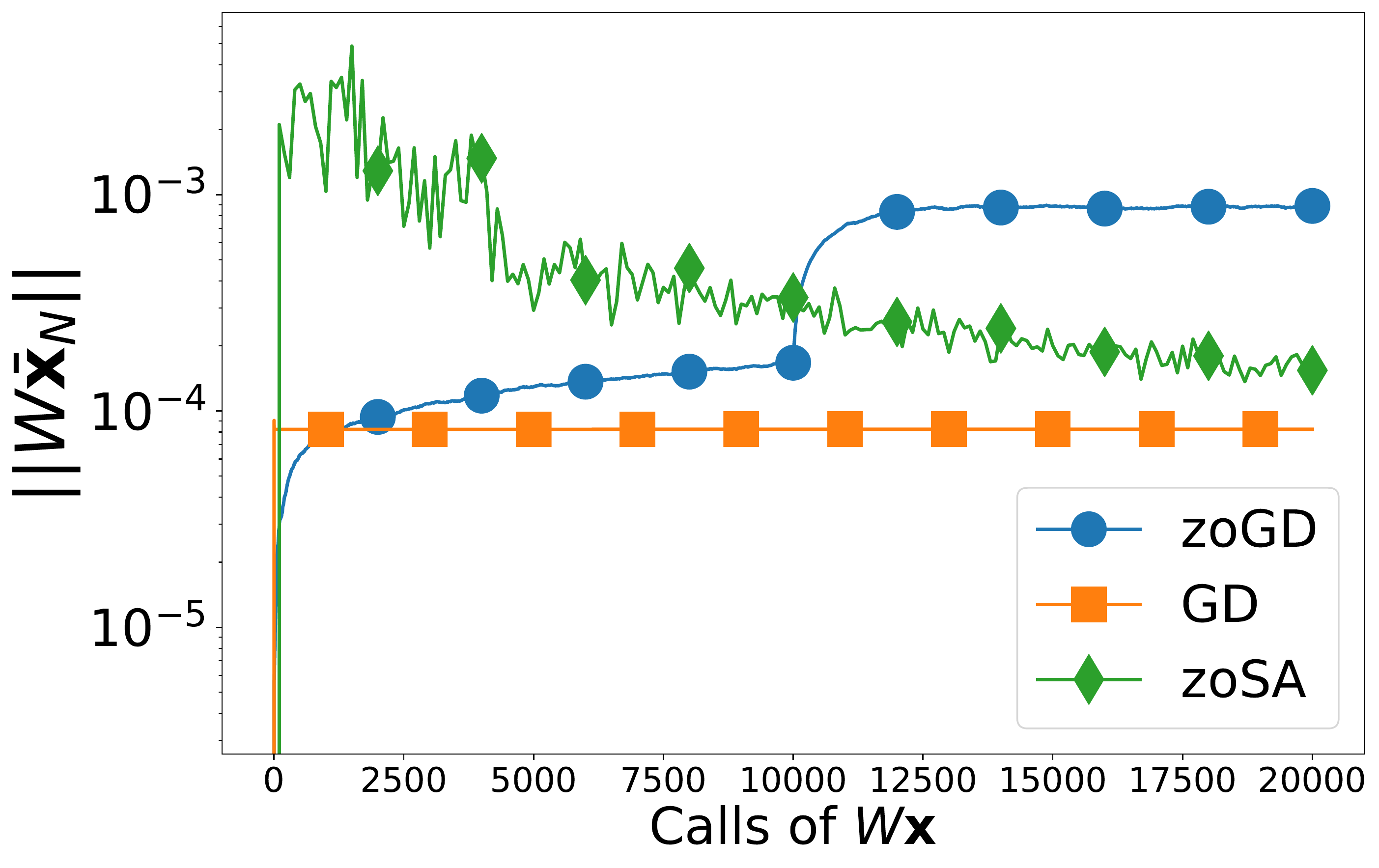}
\includegraphics[width =  0.24\textwidth]{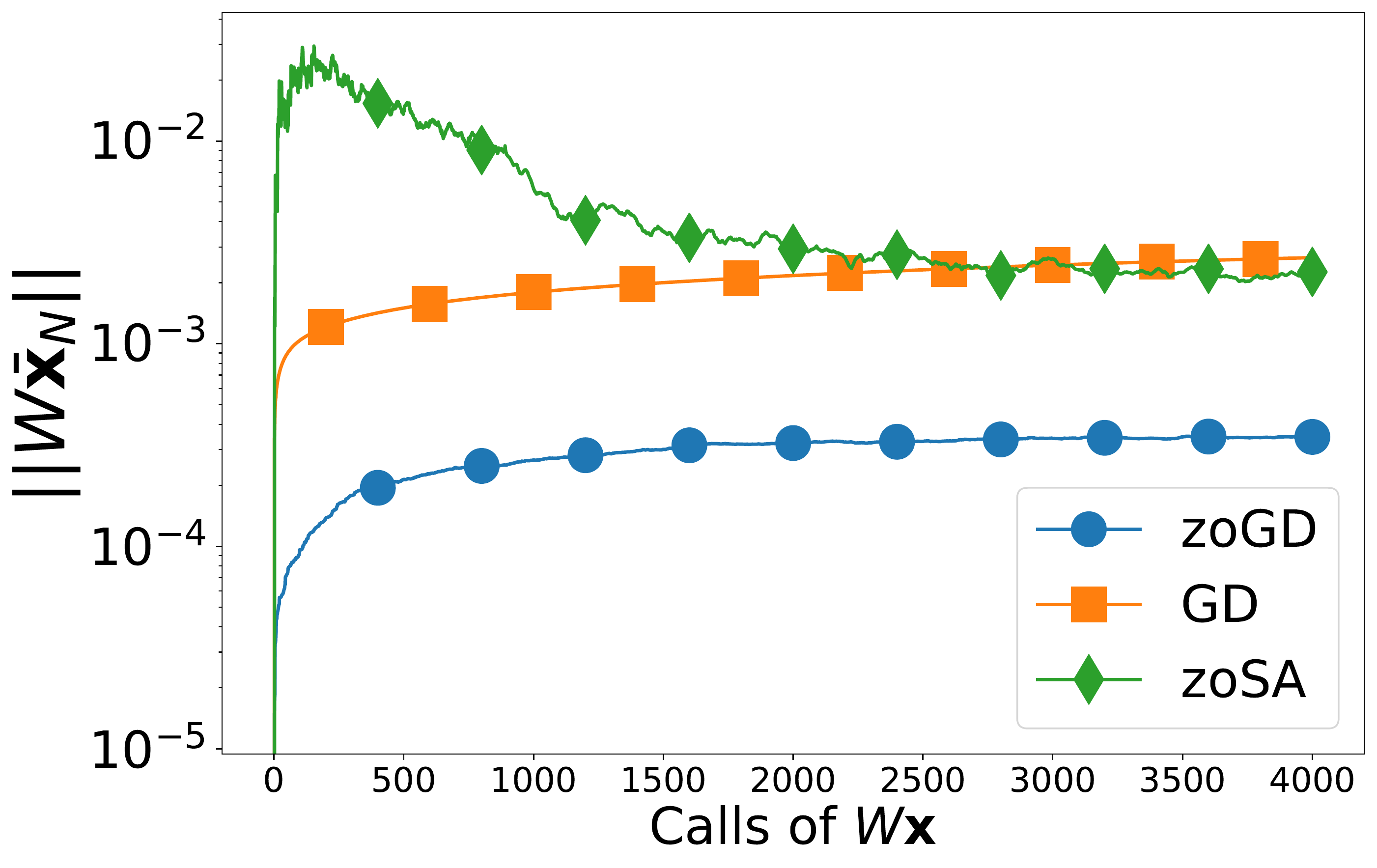}
\includegraphics[width =  0.24\textwidth]{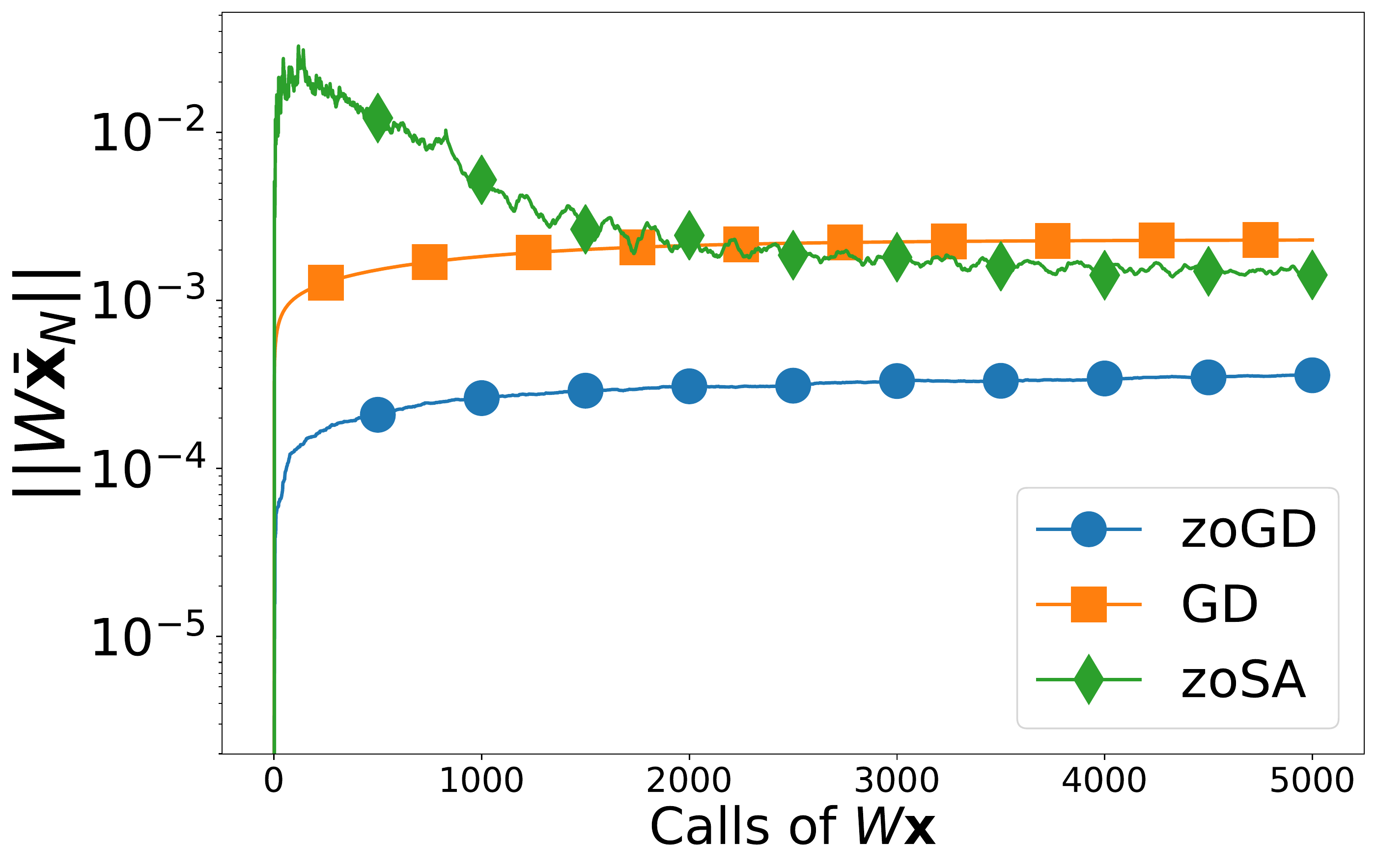}
\begin{minipage}{0.24\textwidth}
\centering
~~~(a) star
\end{minipage}%
\begin{minipage}{0.24\textwidth}
\centering
~~~~~~~(b) complete graph
\end{minipage}%
\begin{minipage}{0.24\textwidth}
\centering
~~~~~~~(c) chain
\end{minipage}%
\begin{minipage}{0.24\textwidth}
\centering
~~~~~~~~~(d) cycle
\end{minipage}%
\caption{{\tt zoSA}, {\tt GD} and {\tt zoGD} applied to solve \eqref{eq:geom_median_problem_distrib} with $R = 10^2$ for different network topologies. First two rows shows how the relative functional gap changes for the methods during their work and the last row shows the evolution of $W\overline{\x}_N$.}
\label{fig:distrib_reg10^2}
\end{figure*}

\begin{figure*}[h!]
\centering
\includegraphics[width =  0.24\textwidth]{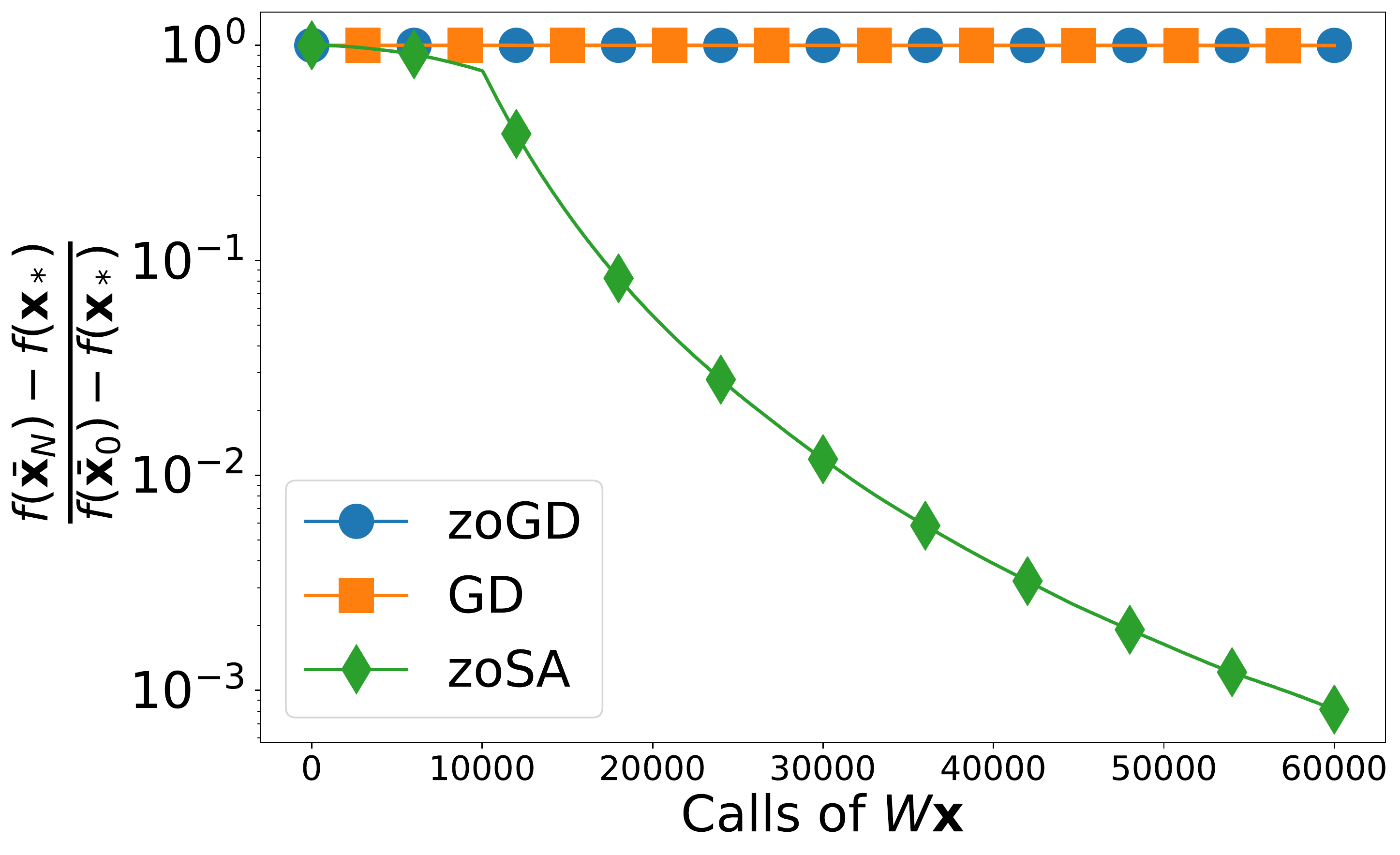}
\includegraphics[width =  0.24\textwidth]{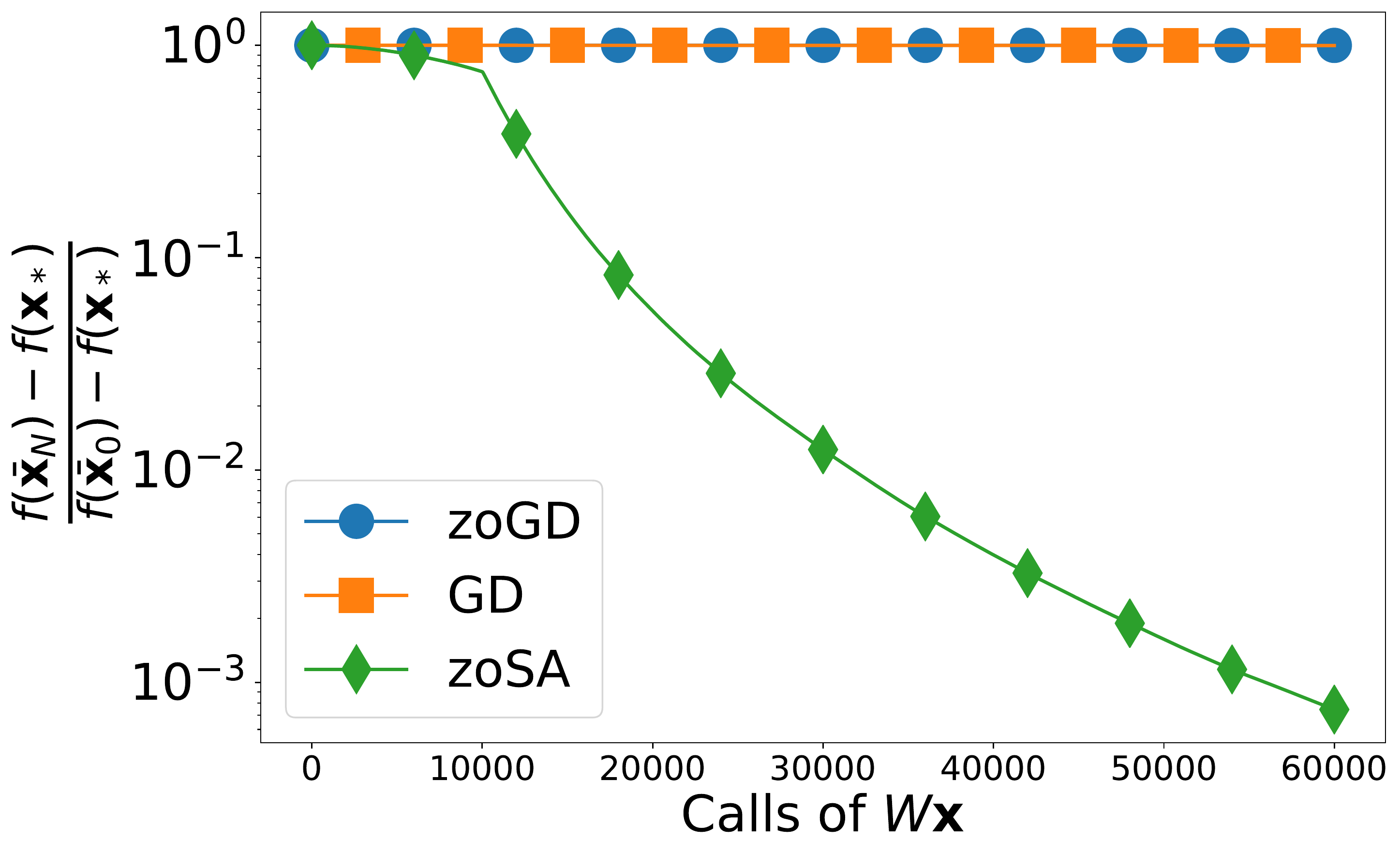}
\includegraphics[width =  0.24\textwidth]{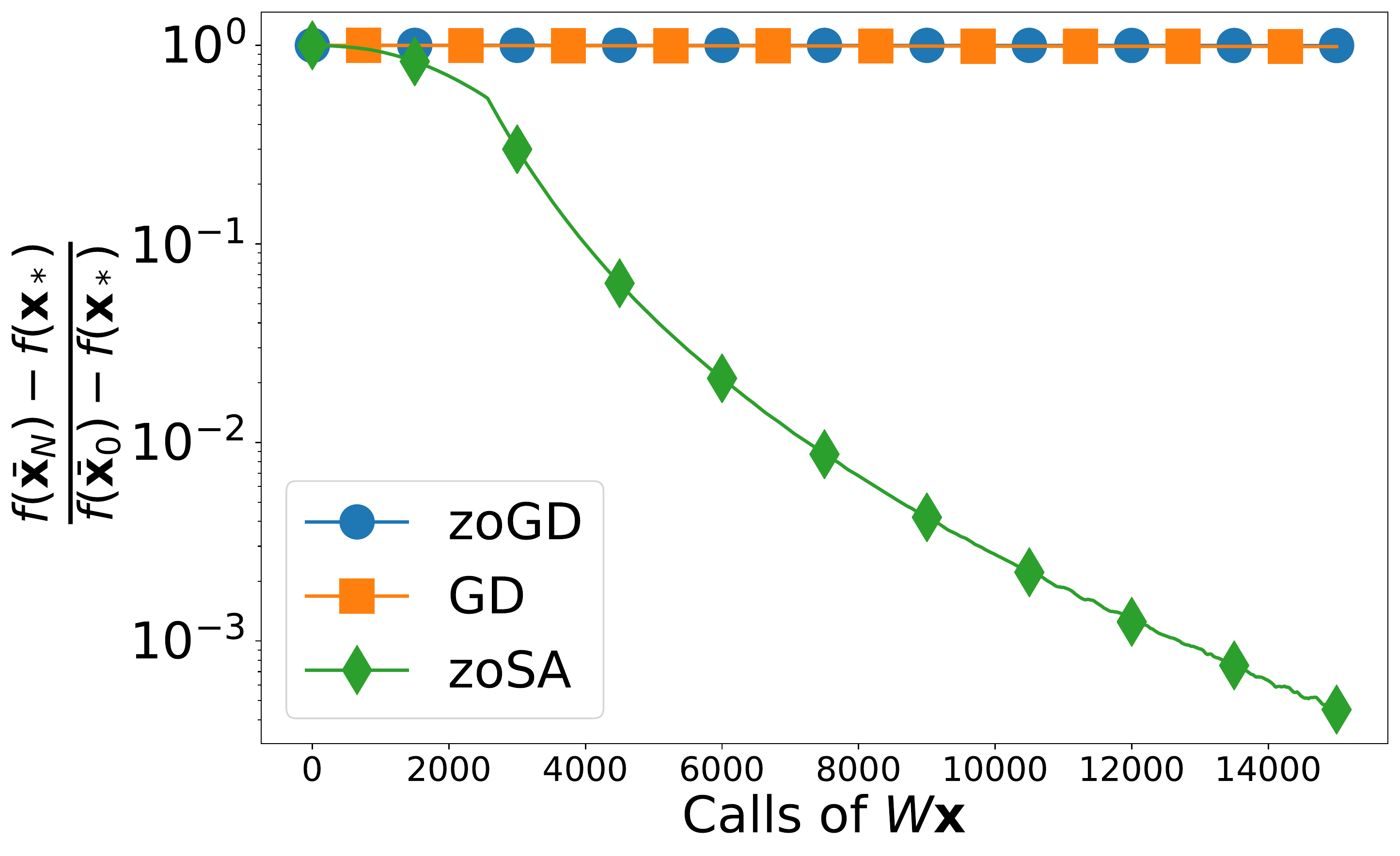}
\includegraphics[width =  0.24\textwidth]{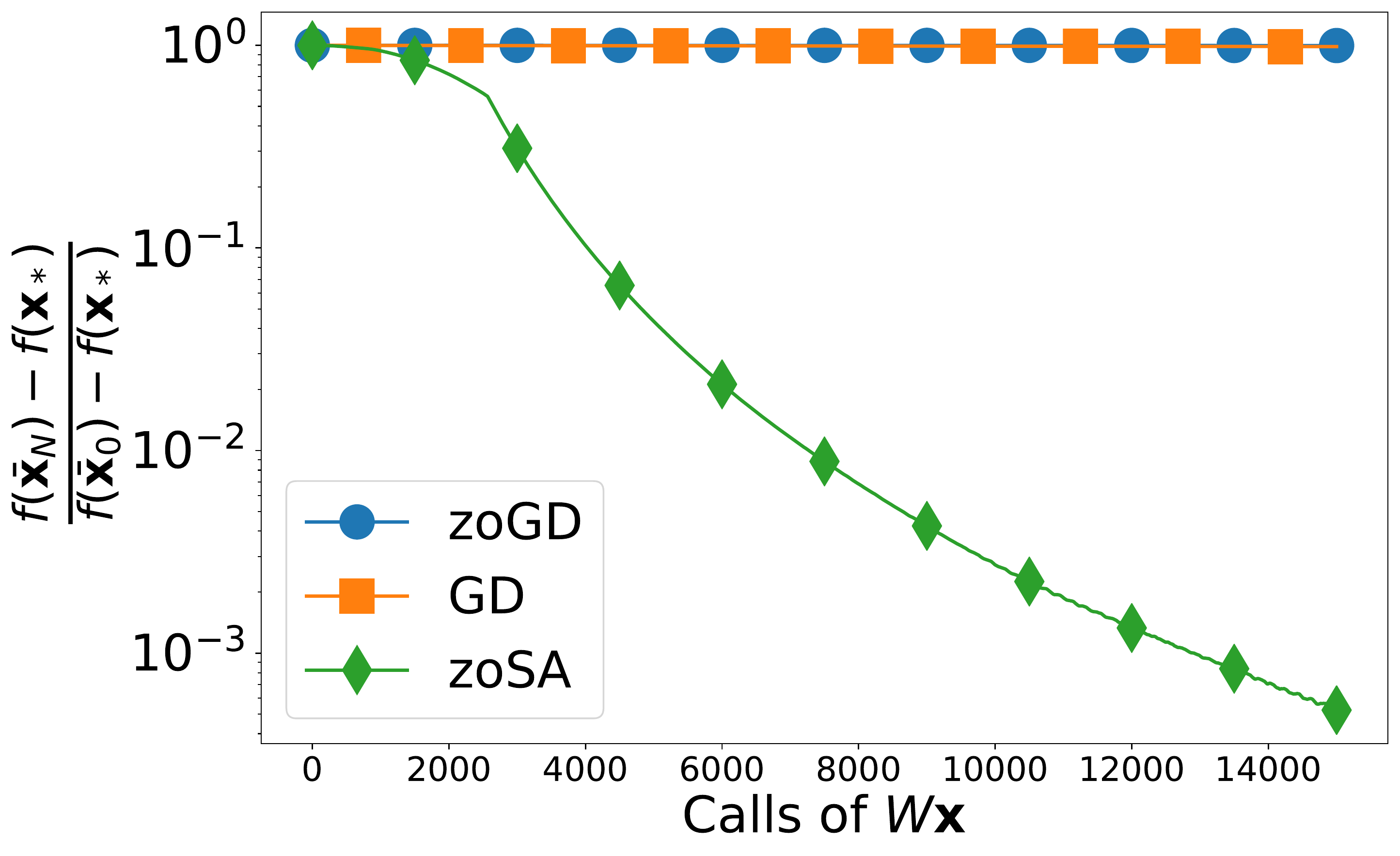}
\includegraphics[width =  0.24\textwidth]{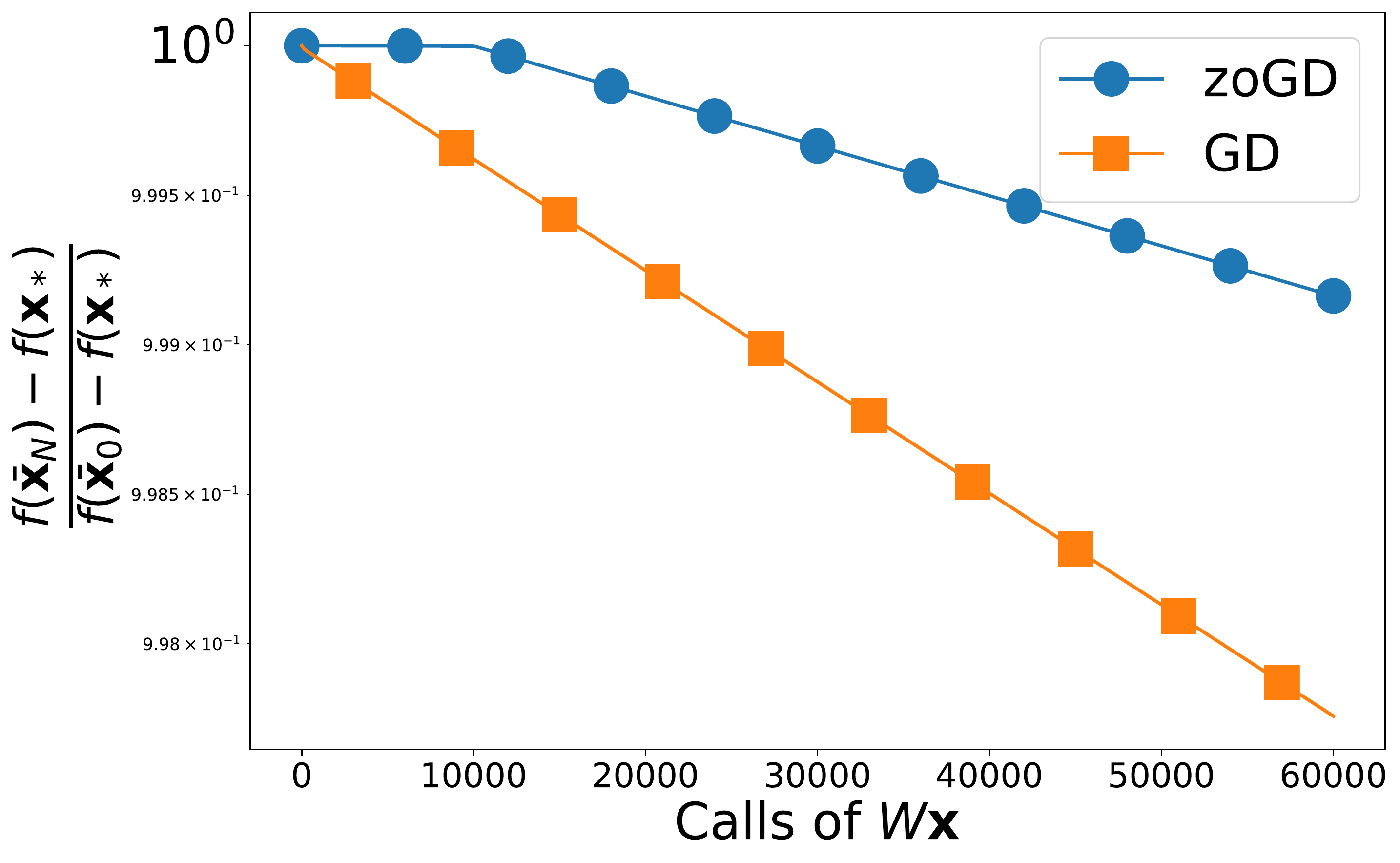}
\includegraphics[width =  0.24\textwidth]{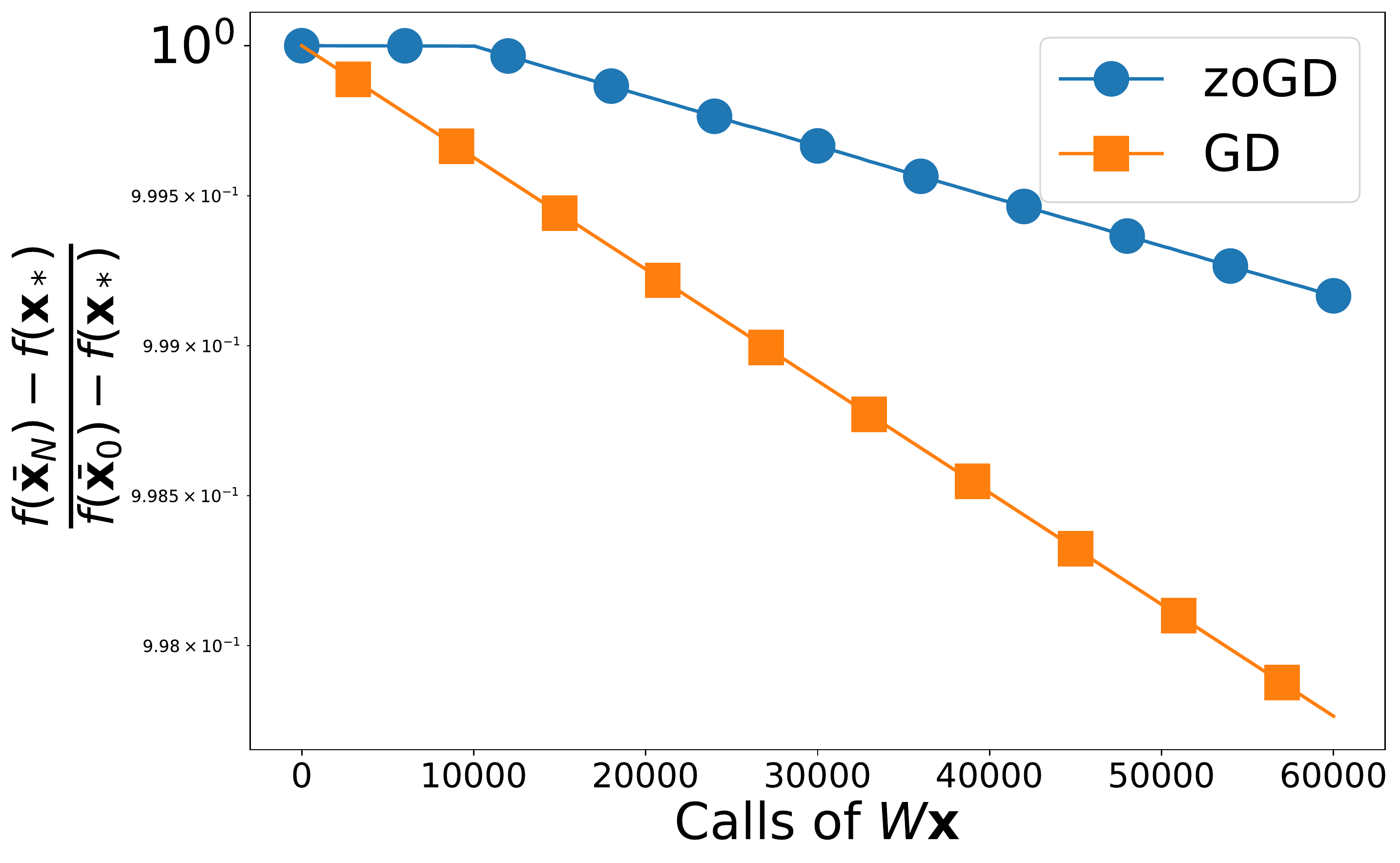}
\includegraphics[width =  0.24\textwidth]{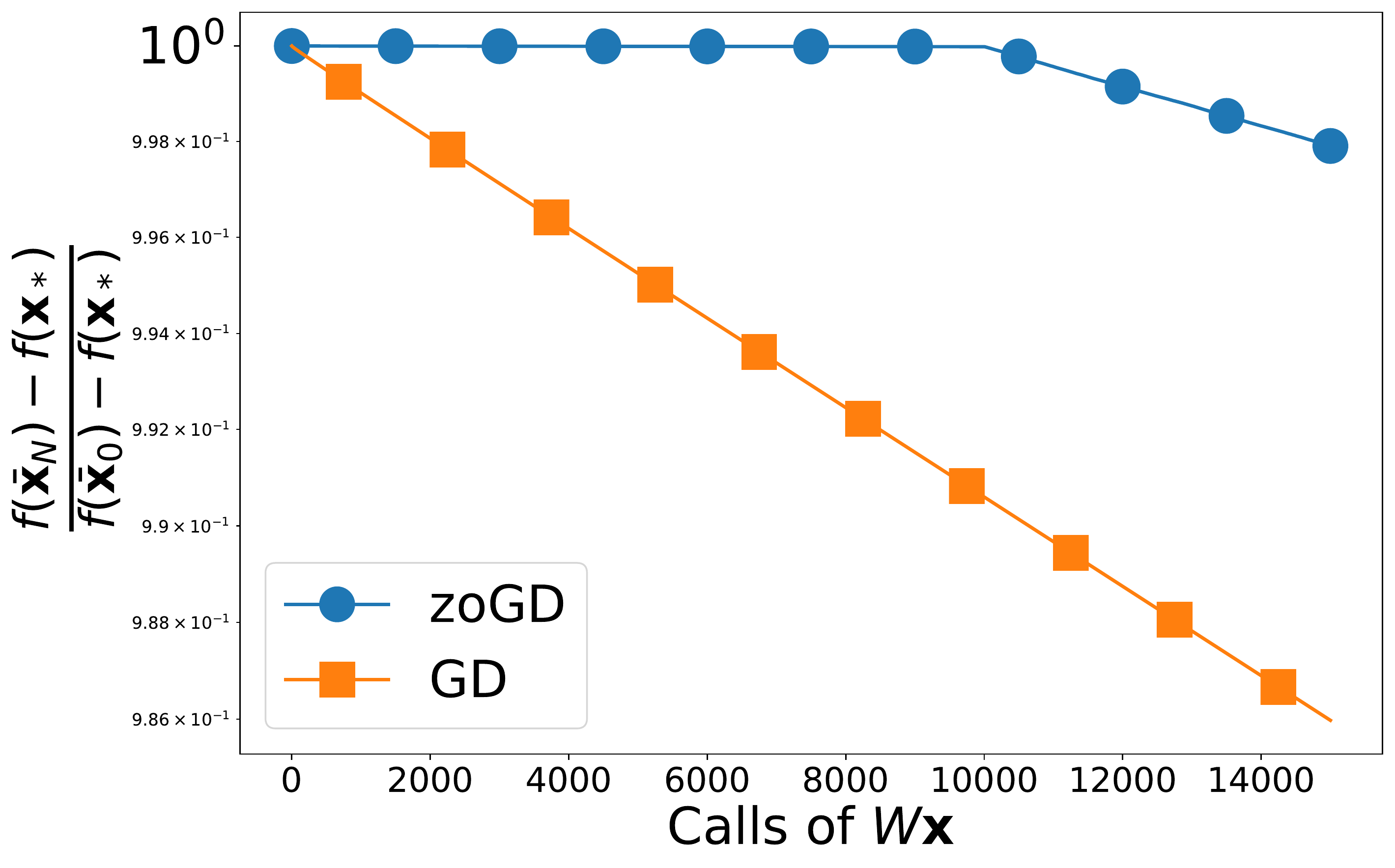}
\includegraphics[width =  0.24\textwidth]{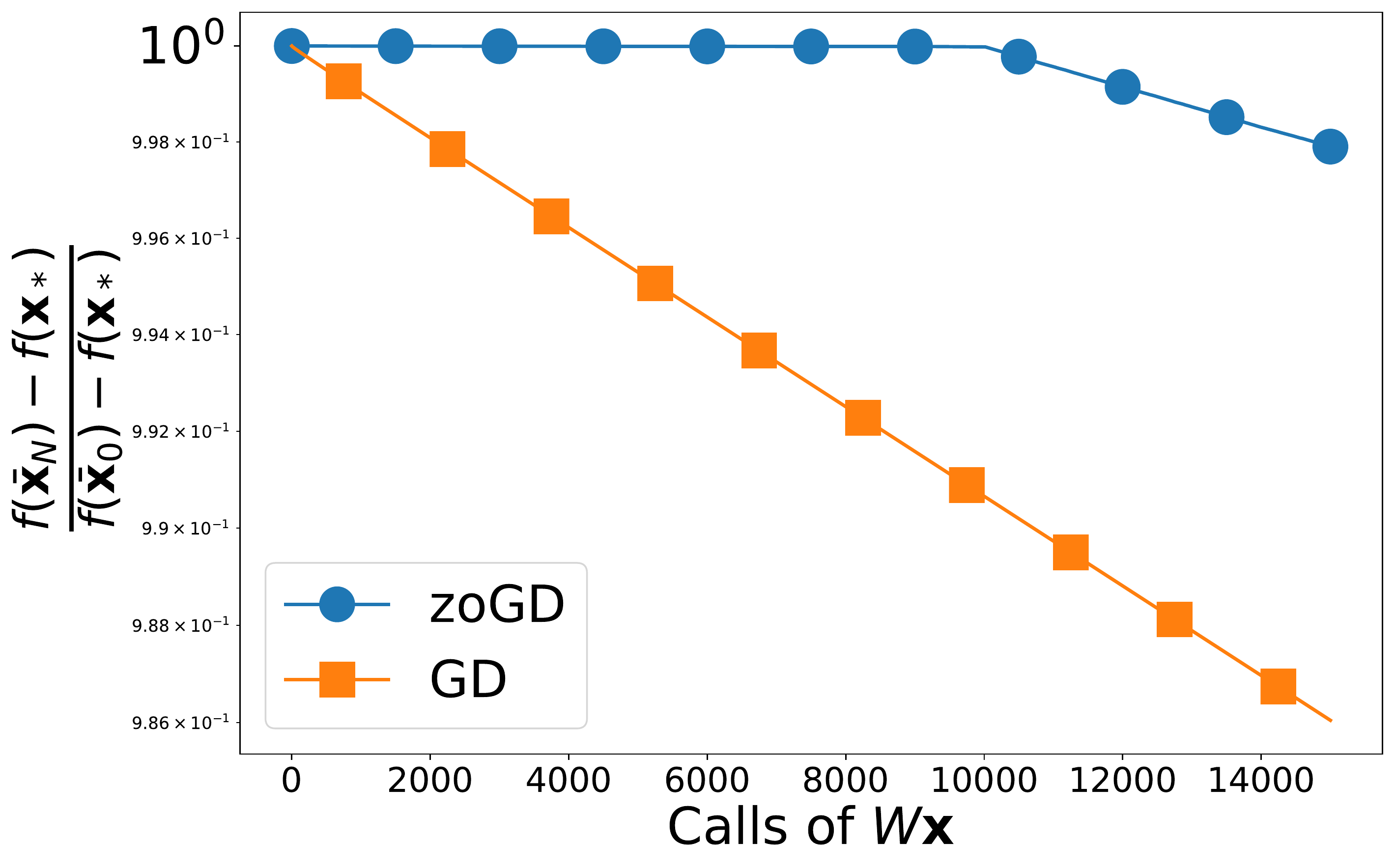}
\includegraphics[width =  0.24\textwidth]{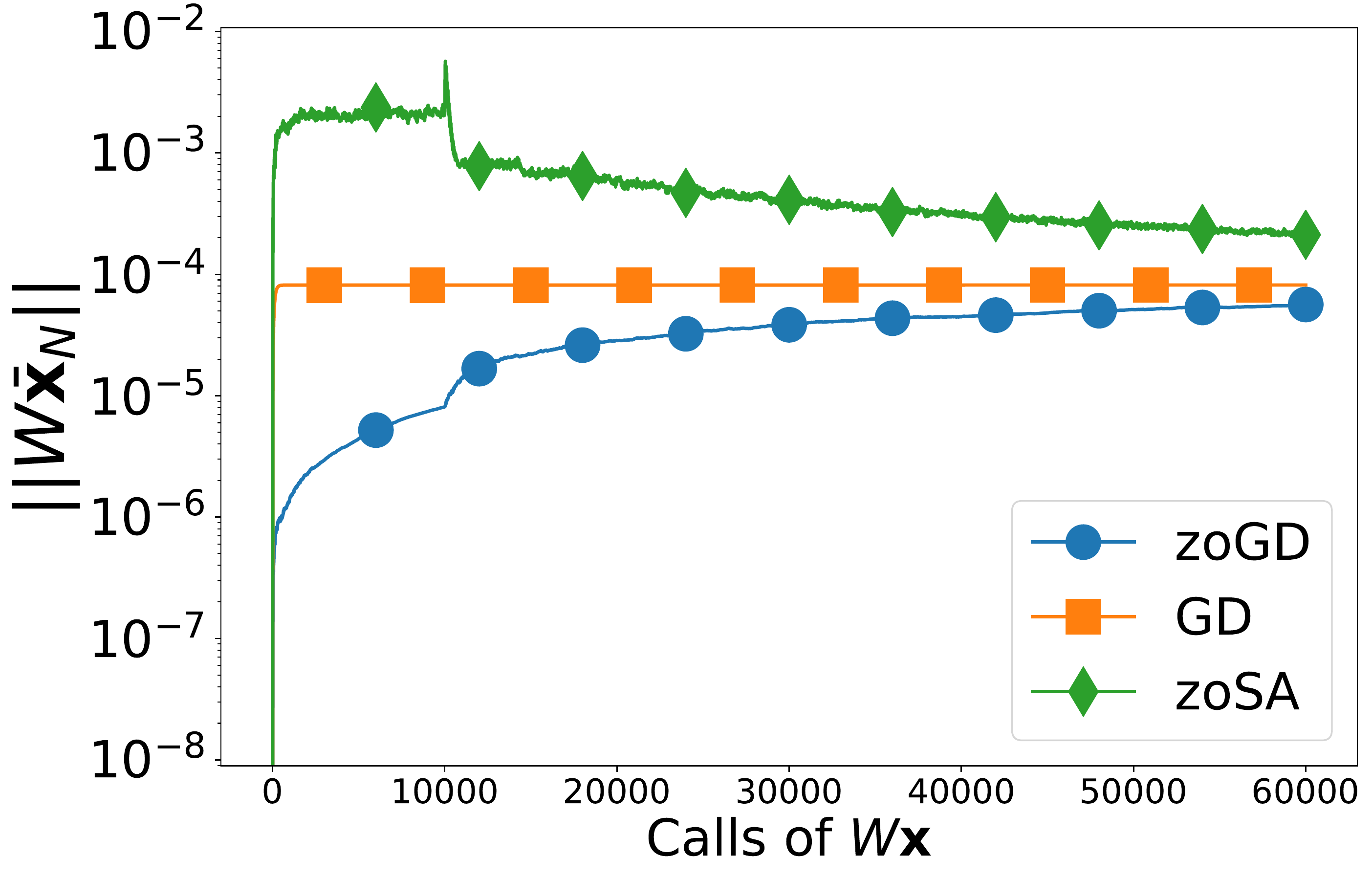}
\includegraphics[width =  0.24\textwidth]{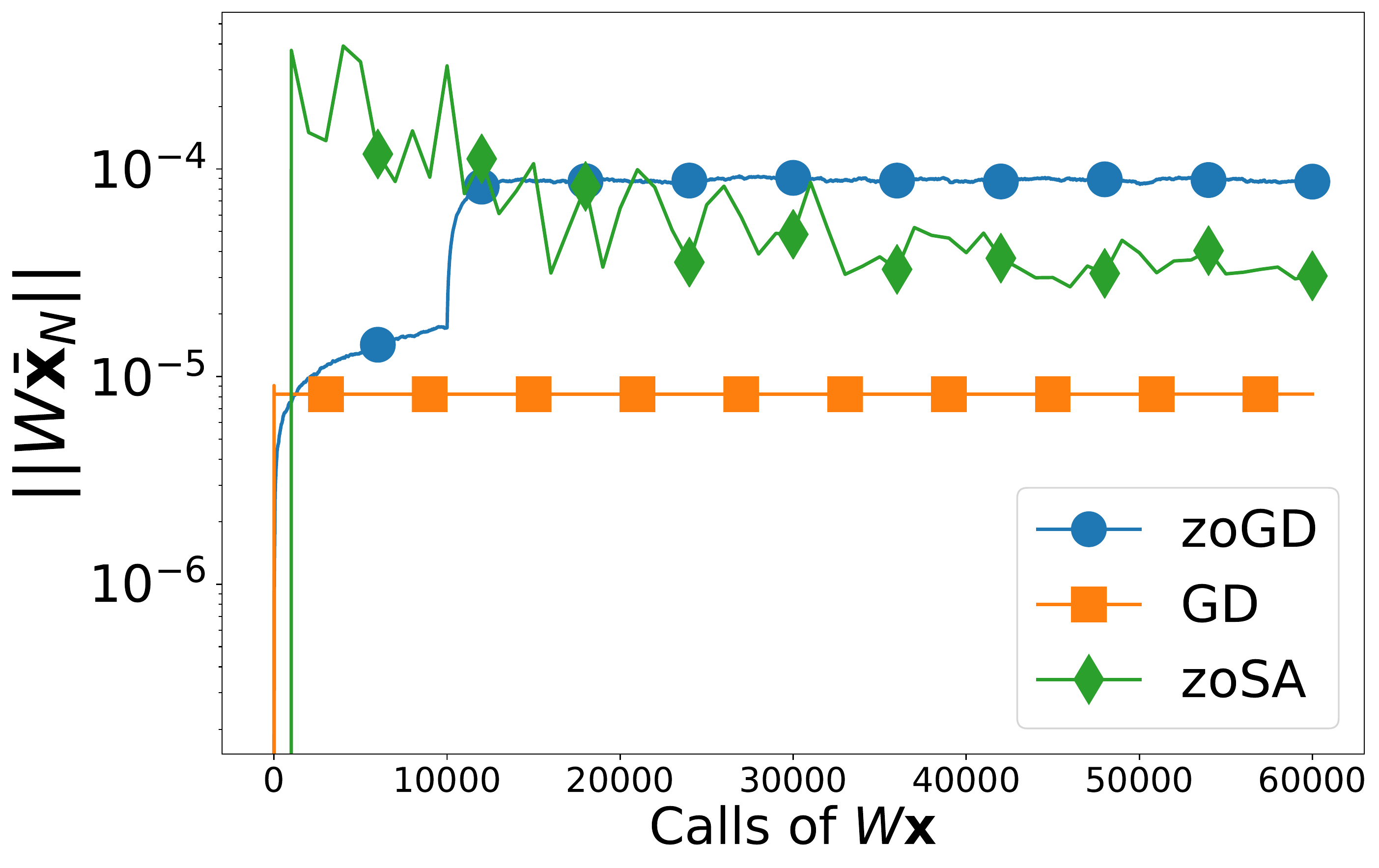}
\includegraphics[width =  0.24\textwidth]{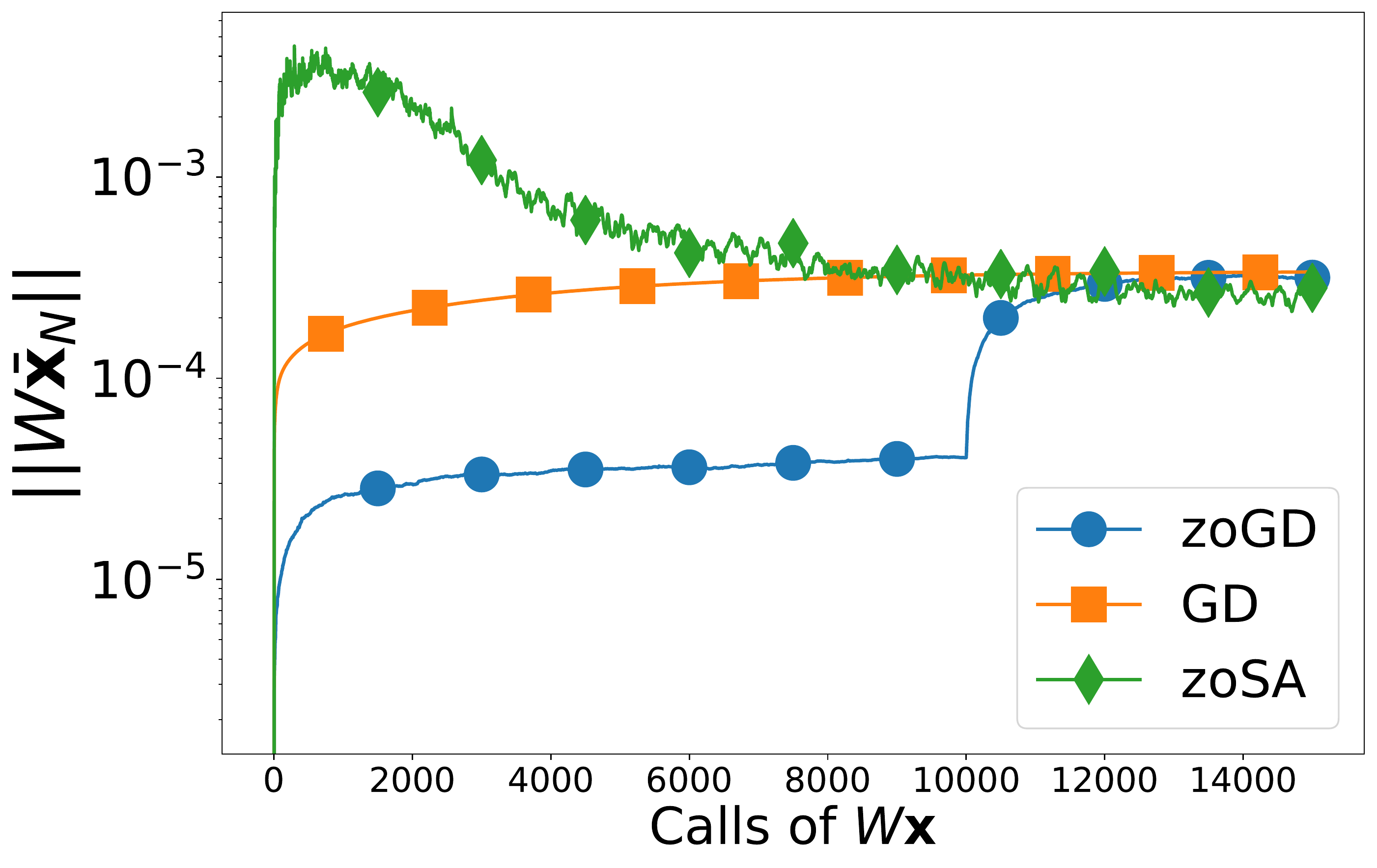}
\includegraphics[width =  0.24\textwidth]{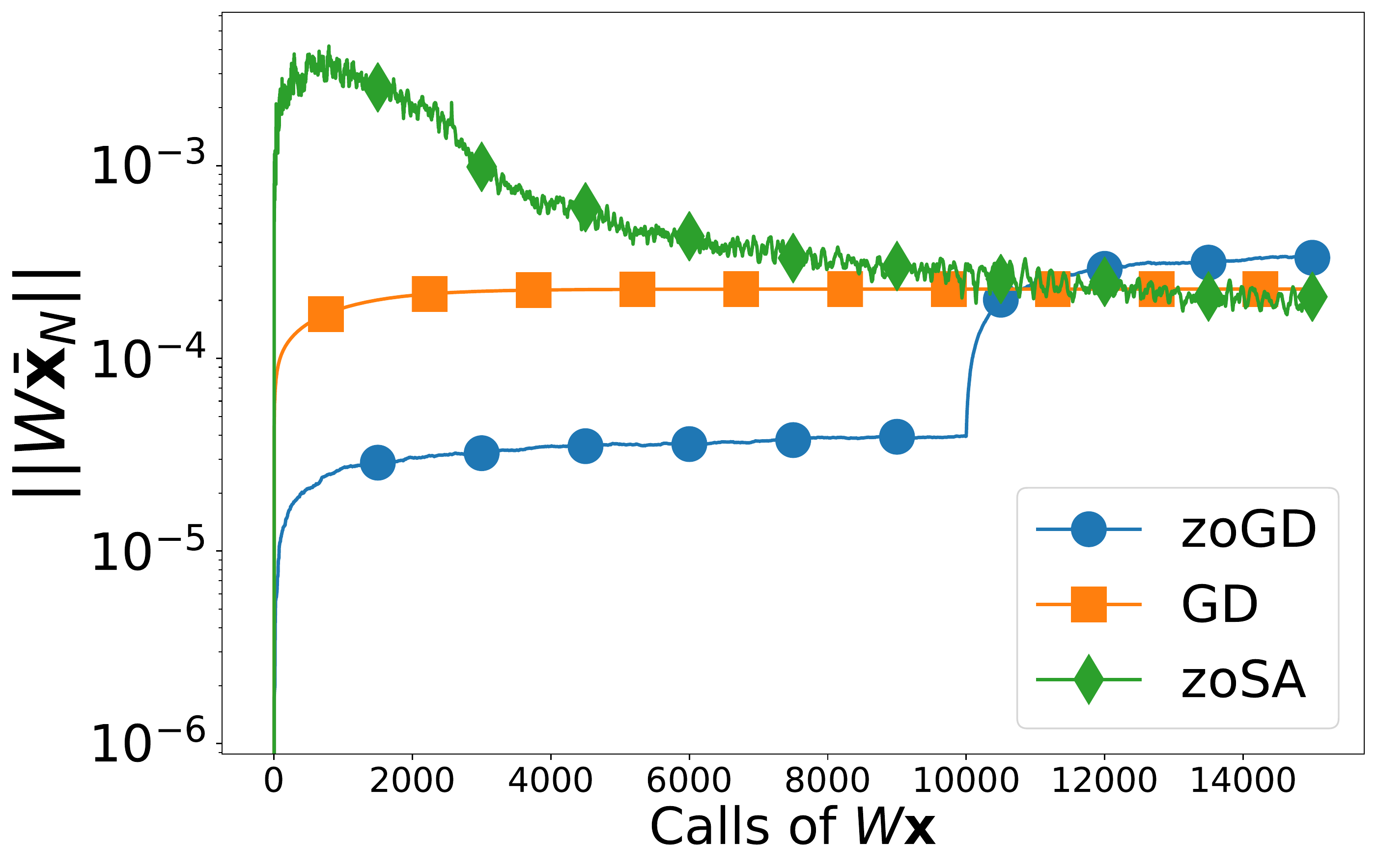}
\begin{minipage}{0.24\textwidth}
\centering
~~~(a) star
\end{minipage}%
\begin{minipage}{0.24\textwidth}
\centering
~~~~~~~(b) complete graph
\end{minipage}%
\begin{minipage}{0.24\textwidth}
\centering
~~~~~~~(c) chain
\end{minipage}%
\begin{minipage}{0.24\textwidth}
\centering
~~~~~~~~~(d) cycle
\end{minipage}%
\caption{{\tt zoSA}, {\tt GD} and {\tt zoGD} applied to solve \eqref{eq:geom_median_problem_distrib} with $R = 10^3$ for different network topologies. First two rows shows how the relative functional gap changes for the methods during their work and the last row shows the evolution of $W\overline{\x}_N$.}
\label{fig:distrib_reg10^3}
\end{figure*}

In our numerical experiments we use a machine with $4$ cores, each is Intel(R) Core(TM) i7-9750H CPU @ 2.60 GHz. We implemented {\tt zoSA}, mirror descent \cite{nemirovsky1983problem} and zeroth-order version of mirror descent \cite{duchi2015optimal} for Euclidean setup, i.e.\ gradient descent ({\tt GD}) and its zeroth-order version ({\tt zoGD}). As we mentioned before, to the best of our knowledge, {\tt zoSA} is the first method for problems of the type \eqref{problem_orig} that uses first-order oracle for the smooth component $g$ and zeroth-order oracle for the non-smooth component $f$. Therefore, we compare {\tt zoSA} with {\tt GD} and {\tt zoGD} that are the state-of-the-art first and zeroth-order methods respectively for convex non-smooth optimization problems.

\subsection{Distributed Computation of Geometric Median}
We consider the problem of searching geometric median \cite{minsker2015geometric, cohen2016geometric} of $m$ vectors $b_1,\ldots,b_m \in \R^n$:
\begin{equation*}
    \min\limits_{x\in\R^n} f(x) = \frac{1}{m}\sum\limits_{i=1}^m \|x - b_i\|_2.
\end{equation*}
Following Section~\ref{sec:affine_and_decentral} we consider the following problem:
\begin{equation}
    \min\limits_{\x\in\R^{nm}} F(\x) = \underbrace{\frac{1}{m}\sum\limits_{i=1}^m \overbrace{\|x_i - b_i\|_2}^{f_i(x_i)}}_{f(\x)} + \underbrace{R\|\sqrt{W}\x\|_2^2}_{g(\x)}.\label{eq:geom_median_problem_distrib}
\end{equation}
As it was mentioned before, if $R = \nicefrac{R_{\y}^2}{\e}$, then $F(\overline\x) - \min_{\x\in\R^{nm}} F(\x) \le \e$ implies $f(\overline\x) - \min_{\sqrt{W}\x = 0}f(\x) \le \e$ and $\|\sqrt{W}\overline \x\|_2 \le \nicefrac{2\e}{R_{\y}}$. However, in practice one can use different choices of $R$ if it offers to get faster such a point $\overline\x$ that $\|\sqrt{W}\overline\x\|_2$ is small enough. In particular, we tried different $R$, but the best result that we obtained are for $R = 10^2$ and $R=10^3$.

In our experiments we emulate the work of the decentralized distributed system with given Laplacian matrix $W$ on one machine in order to demonstrate the performance of {\tt zoSA} on the decentralized distributed optimization problems. That is, we store $\x$ as a long vector and count number of $W\x$ computations since it corresponds to the number of communication rounds in the distributed system. In many real distributed networks communication is a bottleneck, therefore, the number of communication rounds measures, to some degree, the running time of the method.

We run {\tt zoSA}, {\tt GD} and {\tt zoGD} on problem \eqref{eq:geom_median_problem_distrib} with $n=10$ and $m=100$ for several standard topologies like star, cycle, chain, i.e.\ path, and complete graph. We construct vectors $b_1,\ldots, b_m$ as i.i.d.\ samples from normal distributed $\cN(\one, 2I_n)$ with the mean $\one = (1,\ldots,1)^\top$ and the covariance matrix $2I_n$ and use the origin of $\R^{nm}$ as a starting point. One can find the results of our numerical experiments on Figures~\ref{fig:distrib_reg10^2}~and~\ref{fig:distrib_reg10^3}. We notice that in these tests {\tt zoSA} outperforms even {\tt GD} which is a first-order method.

\subsection{Logistic Regression with Lasso Regularization}
\begin{figure}[h]
\centering
\vspace{.3in}
\begin{minipage}{0.3\textwidth}
\includegraphics[width =  \textwidth ]{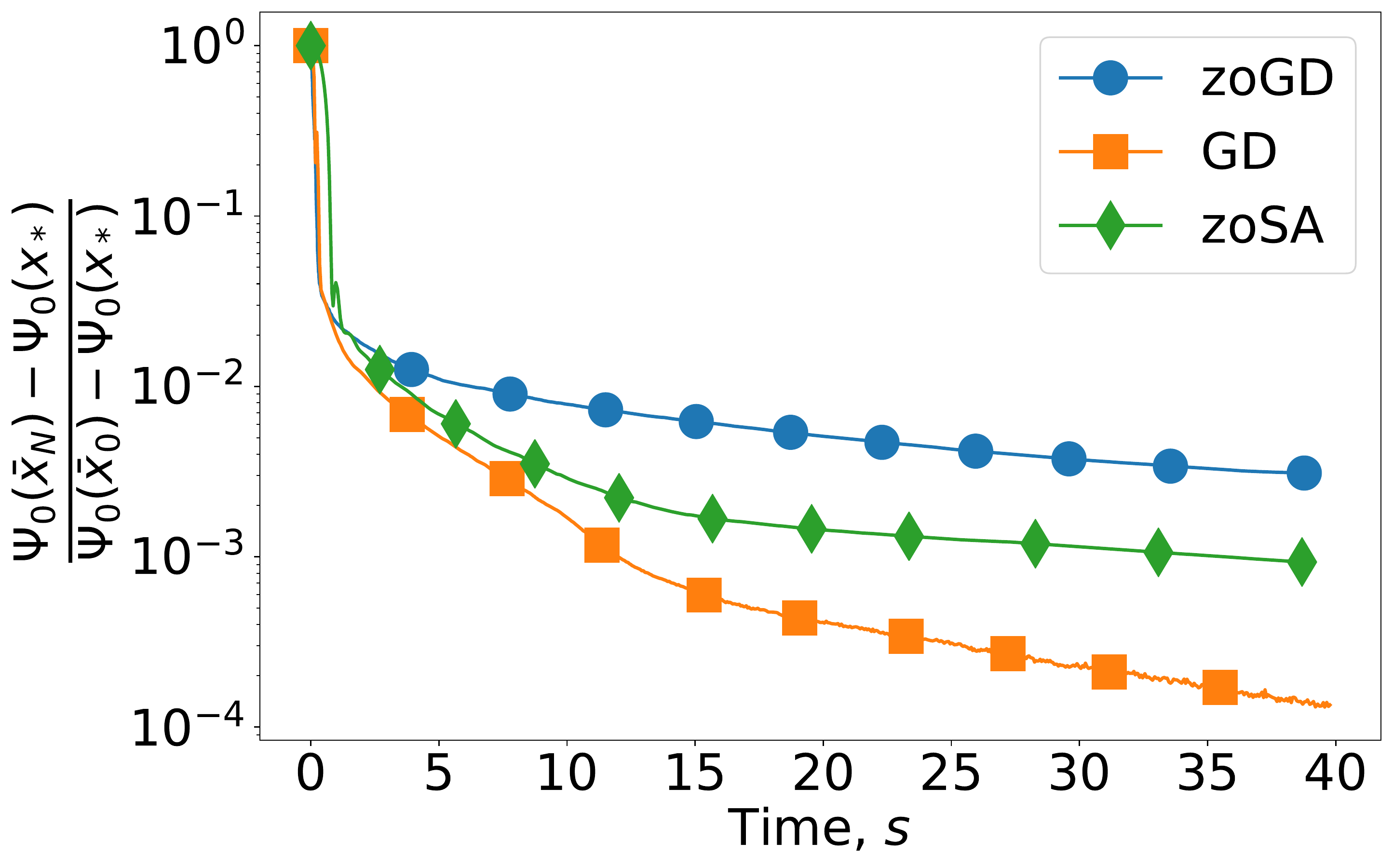}
\end{minipage}%
\begin{minipage}{0.3\textwidth}
  \centering
\includegraphics[width =  \textwidth ]{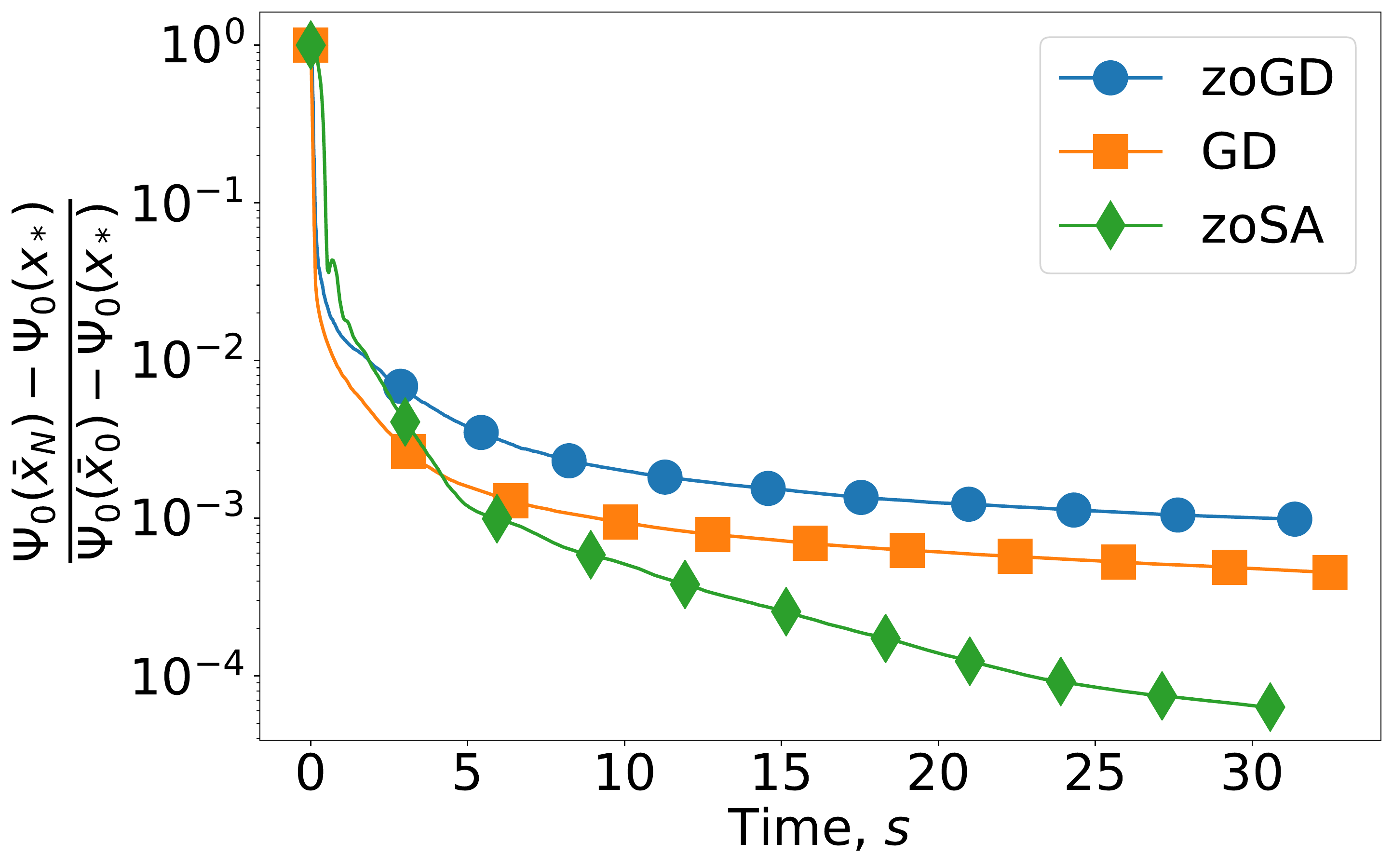}
\end{minipage}%
\begin{minipage}{0.30\textwidth}
  \centering
\includegraphics[width =  \textwidth ]{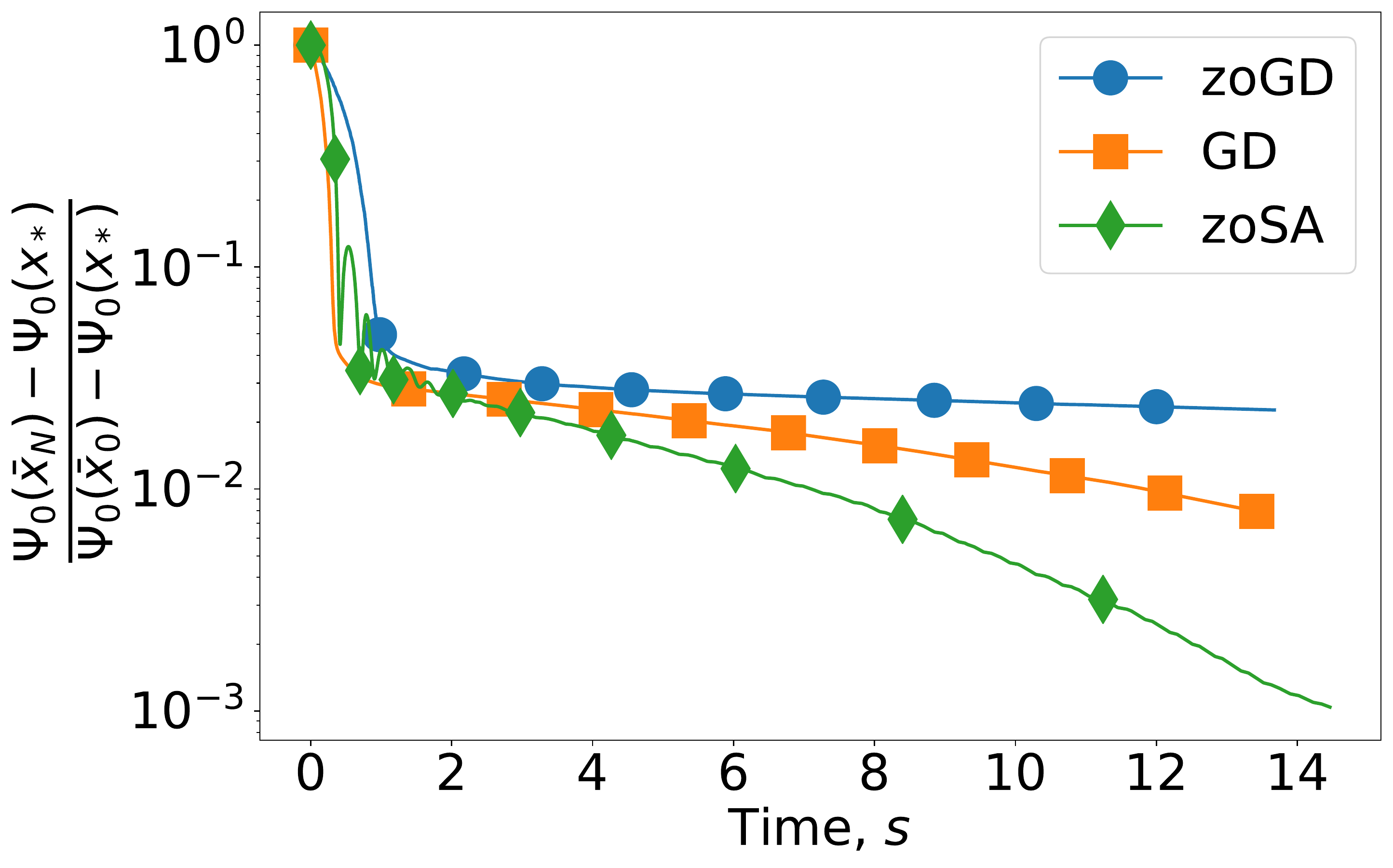}
\end{minipage}%
\\
\begin{minipage}{0.30\textwidth}
\centering
~(a) \texttt{mushrooms}, $l_1 = 10^{-3}$
\end{minipage}%
\begin{minipage}{0.30\textwidth}
  \centering
~~~~~~(b) \texttt{a5a}, $l_1 = 10^{-4}$
\end{minipage}%
\begin{minipage}{0.30\textwidth}
  \centering
(c) \texttt{german.numer}, $l_1 = 10^{-4}$
\end{minipage}%
\caption{{\tt zoSA}, {\tt GD} and {\tt zoGD} applied to solve \eqref{eq:logreg} for different datasets and regularization parameters $l_1$.}
\label{fig:logreg}
\end{figure}

Next, we consider the logistic regression problem with lasso regularization for binary classification:
\begin{eqnarray}
    \min\limits_{x \in \R^n} \Psi_0(x) &=& \overbrace{l_1\|x\|_1}^{f(x)} + g(x) \label{eq:logreg}\\
    g(x) &=& \frac{1}{m}\sum\limits_{i=1}^m \log\left(1 + \exp\left(-y_i\cdot (Ax)_i\right)\right).\notag
\end{eqnarray}
Here $A\in\R^{m\times n}$ is a matrix of objects, $y_1,\ldots,y_m\in\{-1,1\}$ are labels for these objects, $m$ is the size of the dataset and $x\in\R^n$ is a vector of weights. We run {\tt zoSA}, {\tt GD} and {\tt zoGD} on this problem for \texttt{mushrooms} dataset ($m=8124$, $n=112$) with $l_1 = 10^{-3}$, \texttt{a5a} dataset ($m=6414$, $n=123$) with $l_1 = 10^{-4}$ and \texttt{german.numer} dataset ($m=1000$, $n=24$) with $l_1 = 10^{-4}$ \cite{chang2011libsvm}, see Figure~\ref{fig:logreg}. 
For the first case {\tt zoSA} shows the performance that is better than {\tt zoGD}'s performance and worse than {\tt GD}'s one which is reasonable for the method which uses a mixed oracle. However, our method outperforms even {\tt GD} on the second and the third datasets. There is no contradiction here: {\tt zoSA} is based on Sliding Algorithm which has better complexity guarantees than {\tt GD} and {\tt zoSA} has the same complexity as Sliding Algorithm in terms calculations of $\nabla g(x)$.

\subsection{Minimization of Nesterov's Function with Lasso Regularization}
\begin{figure}[h]
    \centering
    \includegraphics[width =  0.6\textwidth ]{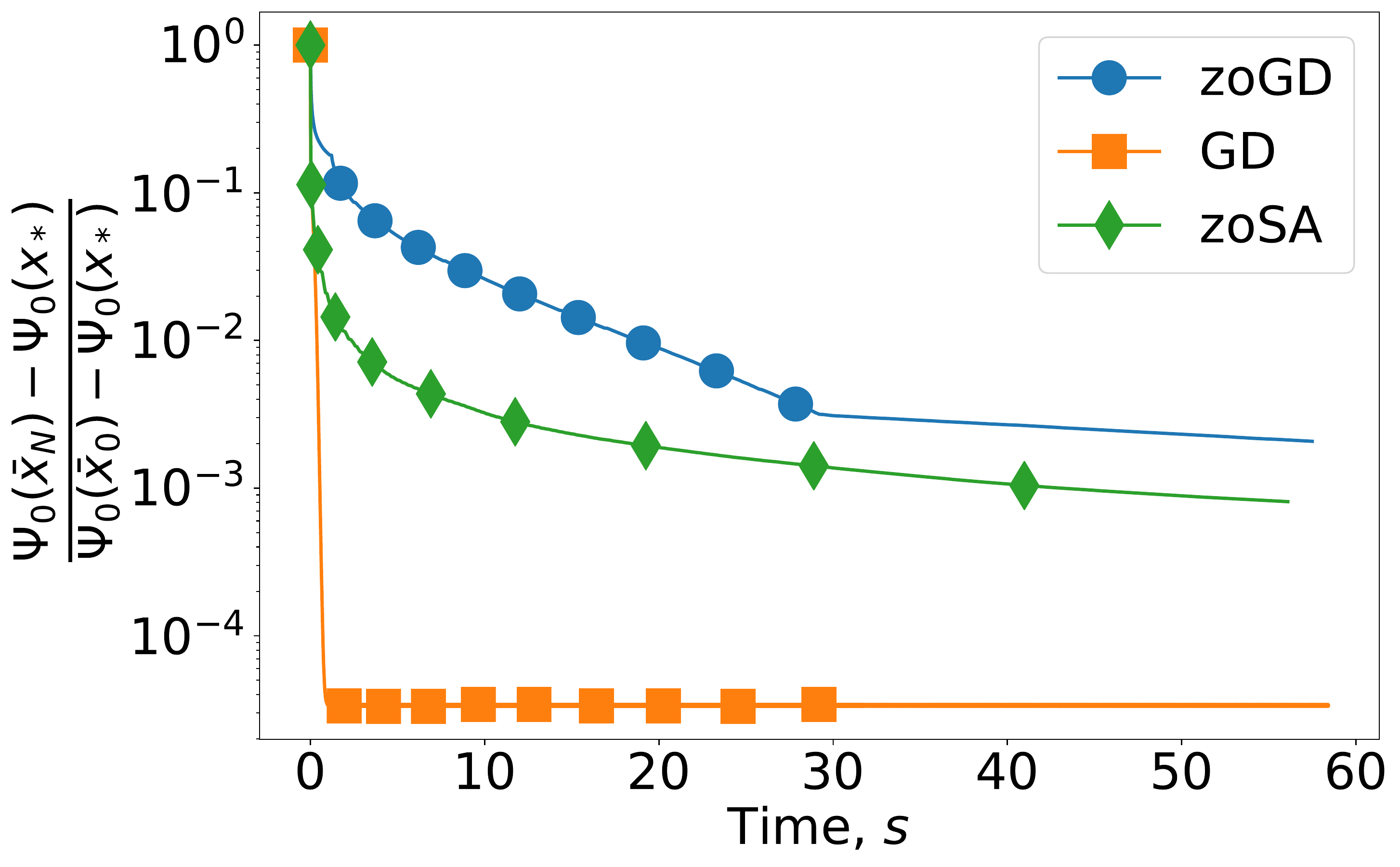}
    \caption{{\tt zoSA}, {\tt GD} and {\tt zoGD} applied to solve \eqref{eq:nesterov_plus_lasso} with $l_1 = 10^{-3}$ and $L = 4$.}
    \label{fig:nesterov_lasso}
\end{figure}

In this section we consider the following problem:
\begin{eqnarray}
    \min\limits_{x \in \R^n} \Psi_0(x) &=& \overbrace{l_1\|x\|_1}^{f(x)} + g(x) \label{eq:nesterov_plus_lasso}\\
    g(x) &=& \frac{L}{8}\left(x_1^2 + \sum\limits_{i=1}^{n-1}(x_i - x_{i+1})^2 + x_n^2\right) - \frac{Lx_1}{4}.\notag
\end{eqnarray}
Here $g(x)$ is a convex and $L$-smooth function, which is one of the ``worst'' functions for the first-order methods in the class of convex and $L$-smooth functions \cite{nesterov2004introduction}, and $f(x)$ has bounded gradients. We run {\tt zoSA}, {\tt GD} and {\tt zoGD} on this problem with $L= 4$ and $l_1 = 10^{-3}$ for a given time. The results are presented in Figure~\ref{fig:nesterov_lasso}. 
Naturally, {\tt zoSA} outperforms {\tt zoGD} since {\tt zoSA} uses first-order oracle for the smooth part while {\tt zoGD} uses only zeroth-order information about $g(x)$. At the same time, our method is inferior to {\tt zoGD} and it is, to some degree, also expected: first-order oracle for $f(x)$ gives more information about descent direction than {\tt zoSA} obtains via zeroth-order oracle.

\section{Basic Facts}
\textbf{Simple upper bound for a squared sum.} For arbitrary integer $n\ge 1$ and arbitrary set of positive numbers $a_1,\ldots,a_n$ we have
\begin{equation}
    \left(\sum\limits_{i=1}^m a_i\right)^2 \le m\sum\limits_{i=1}^m a_i^2\label{eq:squared_sum}
\end{equation}

\textbf{H{\"o}lder inequality.} For arbitrary $x,y\in\R^n$ the following inequality holds
\begin{equation}
    \la x, y \ra \le \|x\|_*\cdot\|y\| \label{eq:holder_ineq}
\end{equation}

\textbf{Cauchy-Schwarz inequality for random variables.} Let $\xi$ and $\eta$ be real valued random variables such that $\EE[\xi^2] < \infty$ and $\EE[\eta^2] < \infty$. Then
\begin{equation}
    \EE[\xi\eta] \le \sqrt{\EE[\xi^2]\EE[\eta^2]}\label{eq:cauchy_schwarz_random}.
\end{equation}

\section{Auxiliary Results}

\begin{lemma}[Lemma 9 from \cite{Shamir15}]\label{lem:lemma_9_shamir} For any function $g$ which  is $L$-Lipschitz with respect to the $\ell_2$-norm, it holds that if $e$ is uniformly distributed on the Euclidean unit sphere, then 
\begin{equation*}
    \sqrt{\mathbb{E}[(g(e) - \mathbb{E}g(e))^4]} \leq c \frac{L^2}{n}
\end{equation*}
for some numerical constant $c$.
\end{lemma}

\begin{lemma}[Lemma 3.5 from \cite{lan}]\label{lem:lemma3.5_lan} Let the convex function $p: X \to \mathbb{R}$, the points $\tilde x, \tilde y \in X$ and scalars $\mu_1, \mu_2 \geq 0$ be given. Let $\nu: X \to \mathbb{R}$ be a differentiable convex function and $V(x,z)$:
$$
V(x,z) = \nu(z) - [\nu(x) + \nabla \nu(x)^\top(z-x)].
$$
If
$$
\tilde u = \argmin_{u \in X}\{p(u) + \mu_1 V(\tilde x, u) + \mu_2 V(\tilde y, u)\}
$$
then for any $u \in X$, we have
\begin{equation*}
p(\tilde u) + \mu_1 V(\tilde x, \tilde u) + \mu_2 V(\tilde y, \tilde u) \leq p(u) + \mu_1 V(\tilde x,  u) + \mu_2 V(\tilde y,  u) - (\mu_1 + \mu_2)V(\tilde u, u).
\end{equation*}
\end{lemma}

\begin{lemma}[Lemma 3.17 from \cite{lan}]\label{lem:lemma3.17_lan}
 Let $w_k \in (0;1]$, $k=1,2,\ldots$ be given. Also let us denote

\begin{align} W_k = 
\begin{cases} 
    1, & k=1,\\
    (1 - w_k)W_{k-1}, & k>1.
\end{cases}\nonumber
\end{align}

Suppose that $W_k>0$ for all $k>1$ and that the sequence $\{\delta_k\}_{k \geq 0}$ satisfies
$$
\delta_k \leq (1 - w_k)\delta_{k-1} + B_k, ~~~ k=1,2,\ldots
$$

for some positive constants $\{B_k\}_{k \geq 0}$.

Then, we have 
\begin{equation*}
\delta_k \leq W_k(1 - w_1)\delta_{0} + W_k \sum\limits_{i=1}^{k} \frac{B_i}{W_i}.
\end{equation*}
\end{lemma}

\section{Missing Proofs from Section~\ref{sec:main_res_cvx}}
\subsection{One Technical Lemma}
\begin{lemma}\label{lem:bounded_grad_to_lan's_condition}
    Assume that for the differentiable function $f$ defined on a closed and convex set $X$ there exists such $M$ that 
    \begin{eqnarray}
    \label{grad_f} 
     \|\nabla f(x) \|_2 \leq M\quad \forall x \in X.
\end{eqnarray}
Then,
\begin{equation*}
    f(x) \le f(y) + \langle \nabla f(y), x - y \rangle + 2MC_1\|x-y\|,\quad \forall x,y \in X.
\end{equation*}
\end{lemma}
\textbf{Proof of Lemma~\ref{lem:bounded_grad_to_lan's_condition}.}
    For arbitrary points $x,y\in X$ we have
    \begin{eqnarray*}
    f(x) &=& f(y) + \int\limits_{0}^{1} \langle \nabla f(y + \tau(x-y)), x-y \rangle d\tau 
     \nonumber\\
     &=& f(y) + \langle \nabla f(y), x-y \rangle + \int\limits_0^1 \langle \nabla f(y + \tau(x-y)) - \nabla f(y), x-y \rangle d\tau \nonumber\\
     &\overset{\eqref{eq:holder_ineq}}{\le}& f(y) + \langle \nabla f(y), x-y \rangle + \int\limits_0^1 \|\nabla f(y + \tau(x-y)) - \nabla f(y) \|_* \cdot \|x-y\| d\tau 
    \nonumber\\
    &\overset{\eqref{eq:norm_equiv},\eqref{grad_f}}{\le}& f(y) + \langle \nabla f(y), x-y \rangle + \int\limits_0^1 2MC_1\|x-y\| d\tau 
    \nonumber\\
    &\leq& f(y) + \langle \nabla f(y), x - y \rangle + 2MC_1\|x-y\|.
\end{eqnarray*}

\subsection{Proof of Lemma~\ref{lem:lemma_8_shamir_main}}
    First of all, Lemma~8 from \cite{Shamir15} implies that $F(x)$ is convex, differentiable and inequality \eqref{grad_sm_main} holds. Next, we use the definition of $F(x)$ and mean value theorem and get that for all $x\in X$
    \begin{eqnarray*}
        |F(x) - f(x)| &=& \left|\EE\left[f(x+re)\right] - f(x)\right| \le \EE\left[|f(x+re) - f(x)|\right] \le \EE\left[\|\nabla f(z(x,x+re))\|_2\cdot\|re\|_2\right] \overset{\eqref{grad_f}}{\le} rM,
    \end{eqnarray*}
    where $z(x,x+re)$ is a convex combination of $x$ and $x+re$. 
    
    Finally, using the symmetry of the distribution of $e$ and \eqref{grad_sm_main} we get
    \begin{eqnarray*}
        \|\nabla F(x)\|_*^2 &=& \left\| \EE\left[\frac{n}{r}f(x+re)e\right]\right\|_*^2 
        \nonumber\\
        &=& \left\| \EE\left[\frac{n}{2r}f(x+re)e\right] +  \EE\left[\frac{n}{2r}f(x-re)\cdot (-e)\right]\right\|_*^2 
        \nonumber\\
        &=& \left\| \EE\left[\frac{n}{2r}\left(f(x+re) - f(x-re)\right)e\right]\right\|_*^2 
        \nonumber\\
        &\leq& \frac{n^2}{4r^2}\EE \left[ \left\|\left(f(x+re) - f(x-re)\right)e\right\|_*^2 \right] 
        \nonumber\\
        &=& \frac{n^2}{4r^2}\mathbb{E}\left[\left(f(x+re)-f(x -re)\right)^2\|e\|_*^2\right]
        \nonumber\\
        &=& \frac{n^2}{4r^2}\mathbb{E}\left[\left(\left(f(x+re)-\alpha\right)-\left(f(x-re)-\alpha\right)\right)^2 \|e\|_*^2\right].
    \end{eqnarray*}
    Next, we apply \eqref{eq:squared_sum} and obtain
    \begin{eqnarray*}
    \frac{n^2}{4r^2}\mathbb{E}\left[\left(\left(f(x+re)-\alpha\right)-\left(f(x-re)-\alpha\right)\right)^2 \|e\|_*^2\right]
    &\leq& 
    \frac{n^2}{2r^2}\mathbb{E}\left[\left(\left(f(x+re)-\alpha\right)^2 +\left(f(x-re)-\alpha\right)^2\right) \|e\|_*^2\right]
    \\&\leq& 
    \frac{n^2}{2r^2}\left(\mathbb{E}\left[ \left(f(x+re)-\alpha\right)^2\|e\|_*^2\right]+\mathbb{E}\left[\left(f(x-re)-\alpha\right)^2 \|e\|_*^2\right]\right).
    \end{eqnarray*}
    
    Since the distribution of $e$ is symmetric
    \begin{eqnarray*}
    \frac{n^2}{2r^2}\left(\mathbb{E}\left[ \left(f(x+re)-\alpha\right)^2 \|e\|_*^2\right]+\mathbb{E}\left[\left(f(x-re)-\alpha\right)^2 \|e\|_*^2\right]\right) 
      &=&\frac{n^2}{r^2}\mathbb{E}\left[\left(f(x+re)-\alpha\right)^2\|e\|_*^2\right] \nonumber\\
      &\overset{\eqref{eq:cauchy_schwarz_random}}{\le}& \frac{n^2}{r^2}\sqrt{\mathbb{E}\left[\|e\|_*^4\right]}\sqrt{\mathbb{E}\left[\left(f(x+re)
  -\alpha\right)^4\right]}.
\end{eqnarray*}
Taking $\alpha=\mathbb{E}[f(x+re)]$ and using Lemma~\ref{lem:lemma_9_shamir} together with the fact that $f(x + re)$ is $Mr$-Lipschitz w.r.t. $e$ in terms of the $\|\cdot\|_2$-norm (since $\|\nabla f(x)\|_2 \le M$) we get 
  \begin{eqnarray*}
      \frac{n^2p_*^2}{r^2}\sqrt{\mathbb{E}\left[\left(f(x+re)
  -\alpha\right)^4\right]} &\leq& \frac{n^2 p_*^2}{r^2} \overline c \frac{(Mr)^2}{n}  = \overline c n p^2_* M^2,
  \end{eqnarray*}
  where $\overline c$ is some positive constant. That is, we proved that
  \begin{eqnarray*}
     \|\nabla F(x)\|_*^2 \leq  \overline c n p_*^2 M^2,
  \end{eqnarray*}
  which implies \eqref{eq:norm_nabla_F_bound_main} with $\tilde{c} = \sqrt{\overline{c}}$.

\subsection{Proof of Lemma~\ref{lem:second_lemma}}
We prove this inequalities in the similar way as it was done in Lemma 10 (see \cite{Shamir15}). Let us start with \eqref{exp_noise_main}:
\begin{eqnarray*}
        \mathbb{E}[\tilde{f}_r'(x)] &=& \frac{n}{2r}\mathbb{E}[(\tilde{f}(x + re) - \tilde{f}(x - re))e] \notag\\
        &=& \frac{n}{2r}\left(\mathbb{E}[f(x + re, \xi) e] - \mathbb{E}[f(x - re, \xi) e] + \mathbb{E}[\Delta(x + re) e] -  \mathbb{E}[\Delta(x - re) e] \right)
\end{eqnarray*}
Taking into account the independence of $e$, $\xi$ and \eqref{xi0} we have $\mathbb{E}[f(x + re, \xi) e] = \mathbb{E}_e\left[\mathbb{E}_\xi[f(x + re, \xi) e]\right] = \mathbb{E}_e\left[f(x + re)]\right]$. Then,
\begin{eqnarray}
\label{temp1}
        \|\mathbb{E}[\tilde{f}_r'(x)] - \nabla F(x)\|_* 
        &=& \left\|\frac{n}{2r}\left(\mathbb{E}_e\left[f(x + re) e\right] - \mathbb{E}_e \left[f(x - re) e\right] + \mathbb{E}[\Delta(x + re) e] -  \mathbb{E}[\Delta(x - re) e] \right) - \nabla F(x)\right\|_*\notag\\
        &\overset{\eqref{grad_sm_main}}{=}&
        \frac{n}{2r}\left\|\mathbb{E}_e[\Delta(x + re)e] -  \mathbb{E}_e[\Delta(x - re)e]\right\|_* \notag\\
        &\overset{\eqref{eq:cauchy_schwarz_random}}{\leq}& \frac{n}{r}\sqrt{\mathbb{E}_e \left[|\Delta(x + re)|^2\right]\cdot \mathbb{E}_e\left[\|e \|^2_*\right]}
\end{eqnarray}
Applying boundedness of $\Delta(x)$ and \eqref{eq:condition_on_u} to \eqref{temp1}, we get \eqref{exp_noise_main}.

Next, we prove the second part of the lemma:
  \begin{eqnarray}
  \mathbb{E}[\|\tilde{f}_r'(x)\|^2_{*}] &=& \mathbb{E}\left[\left\|\frac{n}{2r}\left(\tf(x+r e)-\tf(x-re)\right)e\right\|_*^2\right] \nonumber\\
  &\overset{\eqref{eq:squared_sum}}{\le}& \frac{n^2}{2r^2}\mathbb{E}\left[\|e\|_*^2 \left(f(x+re, \xi)-f(x -re, \xi)\right)^2\right] + \frac{n^2}{2r^2}\mathbb{E}\left[\|e\|_*^2 \left(\Delta(x+re) - \Delta(x-re)\right)^2\right]\notag\\
  &\overset{\eqref{eq:squared_sum}}{\le}& \frac{n^2}{2r^2}\mathbb{E}\left[\|{e}\|_*^2 \left(\left(f(x+re, \xi)-\alpha\right)-\left(f(x-re, \xi)-\alpha\right)\right)^2\right] + \frac{n^2}{r^2}\mathbb{E}\left[\|{e}\|_*^2\left(\Delta^2(x+re) + \Delta^2(x-re)\right)\right]\nonumber\\
  &\overset{\eqref{delta},\eqref{eq:cauchy_schwarz_random}}{\le}& \frac{n^2}{r^2}\mathbb{E}_{e}\left[\|e\|_*^2 \left(\left(f(x+re, \xi)-\alpha\right)^2+\left(f(x-re, \xi)-\alpha\right)^2\right)\right] + \frac{2n^2\Delta^2}{r^2}\sqrt{\mathbb{E}\left[\|{e}\|_*^4\right]}\nonumber\\
  &\overset{\eqref{eq:condition_on_u}}{\le}& \frac{n^2}{r^2}\left(\mathbb{E}\left[\|{e}\|_*^2 \left(f(x+re, \xi)-\alpha\right)^2\right]+\mathbb{E}\left[\|{e}\|_*^2\left(f(x-re, \xi)-\alpha\right)^2\right]\right) + \frac{2n^2p_*^2\Delta^2}{r^2}.\label{eq:technical_second_moment}
  \end{eqnarray}
  Since the distribution of $e$ is symmetric we can rewrite the r.h.s.\ of \eqref{eq:technical_second_moment} in the following way:
  \begin{eqnarray*}
    \frac{n^2}{r^2}\left(\mathbb{E}\left[\|e\|_*^2 \left(f(x+re, \xi)-\alpha\right)^2\right]+\mathbb{E}\left[\|{e}\|_*^2\left(f(x-re, \xi)-\alpha\right)^2\right]\right) +  \frac{2n^2p_*^2\Delta^2 }{r^2}&\nonumber\\
      &\hspace{-5cm}=\frac{2n^2}{r^2}\mathbb{E}\left[\|e\|_*^2\left(f(x+re, \xi)-\alpha\right)^2\right] + \frac{2n^2p_*^2\Delta^2}{r^2}.
\end{eqnarray*}
Taking into account the independence of $e$ and $\xi$ we derive
\begin{eqnarray*}
    \frac{2n^2}{r^2}\mathbb{E}\left[\|e\|_*^2\left(f(x+re, \xi)-\alpha\right)^2\right] + \frac{2n^2p_*^2\Delta^2}{r^2} = \frac{2n^2}{r^2}\mathbb{E}_{\xi}\left[ \mathbb{E}_e \left[\|e\|_*^2\left(f(x+re, \xi)-\alpha\right)^2\right]\right] + \frac{2n^2p_*^2\Delta^2}{r^2}.
\end{eqnarray*}
Next, using Cauchy-Schwarz inequality and \eqref{eq:condition_on_u} we obtain
  \begin{eqnarray*}
  \frac{2n^2}{r^2}\mathbb{E}_{\xi}\left[ \mathbb{E}_e \left[\|e\|_*^2\left(f(x+re, \xi)-\alpha\right)^2\right]\right] + \frac{2n^2p_*^2\Delta^2}{r^2} &\overset{\eqref{eq:cauchy_schwarz_random}}{\le}& \frac{2n^2}{r^2}\mathbb{E}_{\xi}\left[\sqrt{\mathbb{E}_e\left[\|{e}\|_*^4\right]}\sqrt{\mathbb{E}_e\left[\left(f(x+re, \xi)
  -\alpha\right)^4\right] }\right] + \frac{2n^2p_*^2\Delta^2}{r^2} \nonumber \\&\leq& \frac{2n^2p_*^2}{r^2}\mathbb{E}_{\xi}\left[\sqrt{\mathbb{E}_e\left[\left(f(x+re, \xi)
  -\alpha\right)^4\right] }\right] + \frac{2n^2p_*^2\Delta^2}{r^2}.
  \end{eqnarray*}
  In particular, taking $\alpha=\mathbb{E}_e[f(x+re, \xi)]$
  and using Lemma~\ref{lem:lemma_9_shamir} with the fact that $f(x + re, \xi)$ is $rM(\xi)$-Lipschitz w.r.t. $e$ in terms of the $\|\cdot\|_2$-norm we get 
  \begin{eqnarray}
      \frac{2n^2p_*^2}{r^2}\mathbb{E}_{\xi}\left[\sqrt{\mathbb{E}_e\left[\left(f(x+re, \xi)
  -\alpha\right)^4\right] }\right] + \frac{2n^2p_*^2\Delta^2}{r^2} &\leq& \frac{2n^2p_*^2}{r^2}\mathbb{E}_{\xi}\left[ c \frac{r^2 M^2(\xi)}{n}\right] + \frac{2n^2p_*^2\Delta^2}{r^2} 
  \overset{\eqref{xi}}{=} 2p_*^2\left(cnM^2 + \frac{n^2\Delta^2}{r^2}\right)\notag,
  \end{eqnarray}
  where $c$ is some positive constant.


\subsection{Proof of Lemma~\ref{lem:third_lemma}}
 The proof of this lemma completely repeats the proof of Proposition 8.3 of \cite{lan}. However, we put it here for consistency. Consider the following functions: 
\begin{equation*}
l_F(u_{t-1}, u) = F(u_{t-1}) + \langle \nabla F_r(u_{t-1}), u - u_{t-1} \rangle,
\end{equation*}
\begin{equation*}
\tilde l_F(u_{t-1}, u) = F(u_{t-1}) + \langle \tf'_r(u_{t-1}), u - u_{t-1} \rangle.
\end{equation*}
These definitions imply that $\tilde l_F(u_{t-1}, u) - l_F(u_{t-1}, u) = \langle \delta_t, u - u_{t-1} \rangle$ where $\delta_t$ is defined in \eqref{delta_main}. Lemmas~\ref{lem:bounded_grad_to_lan's_condition}~and~\ref{lem:lemma_8_shamir_main} imply $F(u_t) \leq l_F(u_{t-1}, u_t) + \tilde M\| u_t - u_{t-1}\|$, where $\tilde M= c \sqrt{n}C_1M$.  Adding $h(u_t) + \beta V(x, u_t)$ to this inequality and applying \eqref{Phi_main} we obtain
\begin{eqnarray*}
\Phi(u_t) \le h(u_t) + l_F(u_{t-1}, u_t) + \beta V(x, u_t) + \tilde M \|u_t - u_{t-1}\|.\end{eqnarray*}
From $\tilde l_F(u_{t-1}, u) - l_F(u_{t-1}, u) = \langle \delta_t, u - u_{t-1} \rangle$ we have
\begin{eqnarray*}
\Phi(u_t) &\le& h(u_t) + l_F(u_{t-1}, u_t) + \beta V(x, u_t) + \tilde M \|u_t - u_{t-1}\|
\nonumber\\
&=& h(u_t) + \tilde l_F(u_{t-1}, u_t) - \langle \delta_t, u_t- u_{t-1} \rangle 
+ \beta V(x, u_t) + \tilde M \|u_t - u_{t-1}\|\nonumber\\
&\overset{\eqref{eq:holder_ineq}}{\le}& h(u_t) + \tilde l_F(u_{t-1}, u_t) + \beta V(x, u_t) + (\tilde M + \|\delta_t\|_*) \|u_t - u_{t-1}\|.
\end{eqnarray*}

Applying Lemma~\ref{lem:lemma3.5_lan} to \eqref{u_t}, we obtain that for all $u \in X$
\begin{eqnarray*}
h(u_t) + \tilde l_F(u_{t-1}, u_t) + \beta V(x, u_t)+ \beta p_t V(u_{t-1}, u_t) &\le& h(u) + \tilde l_F(u_{t-1}, u) + \beta V(x, u) + \beta p_t V(u_{t-1}, u) - \beta(1 + p_t) V(u_t, u)\nonumber\\
&=& h(u) + l_F(u_{t-1}, u) + \langle \delta_t, u - u_{t-1} \rangle + \beta V(x, u) \nonumber\\
&&\quad+ \beta p_t V(u_{t-1}, u) - \beta(1 + p_t) V(u_t, u)\nonumber\\
&\le& \Phi(u) + \beta p_t V(u_{t-1}, u) - \beta (1 + p_t) V(u_t, u) + \langle \delta_t, u - u_{t-1} \rangle,
\end{eqnarray*}
where the last inequality follows from the convexity of $F$ (see Lemma~\ref{lem:lemma_8_shamir_main}) and \eqref{Phi_main}. Moreover, the strong convexity of $V$ implies that
\begin{eqnarray*}
- \beta p_t V(u_{t-1}, u_t) + (\tilde M+\|\delta_t\|_*) \| u_t - u_{t-1}\| &\le& - \frac{\beta p_t }{2} \|u_t - u_{t-1}\|^2 + (\tilde M + \|\delta_t\|_*) \|u_t - u_{t-1}\|\nonumber\\  
&\le& \frac{\left(\tilde M+\|\delta_t\|_*\right)^2}{2\beta p_t},
\end{eqnarray*}
where the last inequality follows from the simple fact that $- \nicefrac{a t^2}{2} + b t \le \nicefrac{b^2}{(2a)}$ for any $a > 0$.

Combining previous three inequalities, we conclude that
\begin{equation*}
\Phi(u_t)  - \Phi(u)  \leq  \beta p_t V(u_{t-1}, u) - \beta (1 + p_t) V(u_t, u) 
+ \frac{\left(\tilde M + \|\delta_t\|_*\right)^2}{2\beta p_t}
+ \langle \delta_t, u - u_{t-1} \rangle.
\end{equation*}
Now dividing both sides of the above inequality by $1 + p_t$ and rearranging the terms, we obtain
\begin{equation*}
\beta V(u_t, u) + \frac{\Phi(u_t)  - \Phi(u)} {1+p_t} 
\leq \frac{\beta p_t}{1+ p_t} V(u_{t-1}, u) + \frac{\left(\tilde M + \|\delta_t\|_*\right)^2}{2 \beta (1+p_t) p_t }
+ \frac{\langle \delta_t, u - u_{t-1} \rangle}{1 + p_t},
\end{equation*}
which, in view of Lemma~\ref{lem:lemma3.17_lan}, implies that
\begin{eqnarray} 
\frac{\beta}{P_t} V(u_t, u) + \sum_{i=1}^t \frac{\Phi(u_i)  - \Phi(u)} {P_i (1+p_i)}
\leq \beta V(u_{0}, u) +  \sum_{i=1}^t \left[ \frac{\left(\tilde M + \|\delta_i\|_*\right)^2}{2 \beta P_i (1+p_i) p_i } 
+ \frac{\langle \delta_i, u - u_{i-1} \rangle}{P_i (1 + p_i)} \right]. \label{temp22}
\end{eqnarray}
By definition of $\tilde u_t$ (see \eqref{tilde_u_t}) and \eqref{p_t_main} we have
\begin{eqnarray}
\label{temp34}
\tilde u_t &=& \frac{P_t}{1 - P_t} \left(\frac{1 - P_{t-1}}{P_{t-1}}\tilde u_{t-1} + \frac{1}{P_t(1 + p_t)}u_t \right),\nonumber\\
\tilde u_t &=& \frac{P_t}{1 - P_t} \left(\frac{1 -P_{t-2}}{P_{t-2}}\tilde u_{t-2} + \frac{1}{P_{t-1}(1 + p_{t-1})}u_{t-1} + \frac{1}{P_t(1 + p_t)}u_t \right) = \ldots = \frac{P_{t}}{1 - P_{t}}\sum\limits_{i=1}^t \frac{1}{P_i(1 + p_i)}u_i.
\end{eqnarray}
Combining \eqref{temp22} and \eqref{temp34} we get the result.

\subsection{Proof of Theorem~\ref{thm:first_theorem}}
    The proof of this theorem is almost identical to the proof of Theorem 8.2 from \cite{lan} and via performing similar steps one can get the following inequality which is an analogue of inequality (8.1.69) from \cite{lan}. For convenience, we put below the full proof.

    Using \eqref{lemma_2_main}, definition of $\Phi_k$ and $(x_k, \tilde x_k)$ we have that for all $u\in X$
\begin{eqnarray}
        \label{1.65}
        \beta_k (1 - P_{T_k})^{-1}V(x_k, u) + \left[\Phi_k(\tilde x_k) - \Phi_k(u)\right] &\leq& \beta_k P_{T_k}(1 - P_{T_k})^{-1}V(x_{k-1},u)
        \nonumber\\
        &&\quad + \frac{P_{T_k}}{1 - P_{T_k}} \sum\limits_{i=1}^{T_k}  \frac{\frac{(\tilde M + \| \delta_{k,i}\|_*)^2}{2 \beta_k p_i} + \langle \delta_{k,i}, u-u_{k,i-1} \rangle}{p_i P_{i-1}}
\end{eqnarray}
    First, notice that by the definition of $\overline x_k$ and $\underline{x}_k$, we have $\overline x_k - \underline{x}_k = \gamma_k (\tilde x_k - x_{k-1})$. Using this observation, $L$-smoothness of $g$ (see \eqref{g-L-smooth}), the definition of $l_g$ in \eqref{lg} and the convexity of $g$, we obtain
\begin{eqnarray*}
g(\overline x_k) &\leq& l_g(\underline{x}_k, \overline x_k) + \frac{L}{2} \|\overline x_k - \underline{x}_k\|^2 \nonumber \\
&=& (1-\gamma_k) l_g(\underline{x}_k, \overline x_{k-1}) + \gamma_k l_g(\underline{x}_k, \tilde x_k) + \frac{L \gamma_k^2}{2} \|\tilde x_k - x_{k-1}\|^2\nonumber \\
&\leq& (1-\gamma_k) g(\overline x_{k-1}) + \gamma_k \left[ l_g(\underline{x}_k, \tilde x_k) +\beta_k V(x_{k-1}, \tilde x_k) \right] \nonumber\\
&&\quad- \gamma_k \beta_k V(x_{k-1}, \tilde x_k) + \frac{L \gamma_k^2}{2} \|\tilde x_k - x_{k-1}\|^2\nonumber \\
&\leq& (1-\gamma_k) g(\overline x_{k-1}) + \gamma_k \left[ l_g(\underline{x}_k, \tilde x_k) +\beta_k V(x_{k-1}, \tilde x_k) \right] \nonumber\\
&&\quad - \left( \gamma_k \beta_k - L \gamma_k^2 \right) V(x_{k-1}, \tilde x_k) \nonumber \\
&\leq& (1-\gamma_k) g(\overline x_{k-1}) + \gamma_k \left[ l_g(\underline{x}_k, \tilde x_k) +\beta_k V(x_{k-1}, \tilde x_k) \right],
\end{eqnarray*}
where the third inequality follows from the strong convexity of $V$ and the last inequality follows from \eqref{gamma_kk_main}. By the convexity of $F$, we have 
\begin{equation*}
F(\overline x_k) \leq (1-\gamma_k) F(\overline x_{k-1})  + \gamma_k F(\tilde x_k).
\end{equation*}

Summing up previous two inequalities, and using the definitions of $\Psi$ and $\Phi_k(\tilde{x}_k) \eqdef F(\tilde{x}_k) + l_g(\underline{x}_k, \tilde x_k) +\beta_k V(x_{k-1}, \tilde x_k)$, we have
\[
\Psi(\overline x_k) 
\leq (1-\gamma_k) \Psi(\overline x_{k-1}) + \gamma_k \Phi_k(\tilde x_k).
\]
Subtracting $\Psi(u)$ from both sides of the above inequality,
we obtain
\begin{equation*}
\Psi(\overline x_k) - \Psi(u) \leq (1- \gamma_k) [\Psi(\overline x_{k-1})  - \Psi(u)] + \gamma_k [\Phi_k(\tilde x_k) - \Psi(u)].
\end{equation*}
Also note that by the definition of $\Phi_k$ and the convexity of $g$,
\begin{equation*}
\Phi_k(u) \le F(u) + g(u) + \beta_k V(x_{k-1},u) = \Psi(u) + \beta_k V(x_{k-1}, u), \ \ \forall u \in X.
\end{equation*}
Combining these two inequalities, we obtain for all $u \in X$
\begin{eqnarray}
\Psi(\overline x_k) - \Psi(u) \leq (1- \gamma_k) [\Psi(\overline x_{k-1})  - \Psi(u)] + \gamma_k [\Phi_k(\tilde x_k) - \Phi_k(u) +  \beta_k V(x_{k-1}, u)]. \label{bnd_outer1}
\end{eqnarray}

Using \eqref{1.65} and \eqref{bnd_outer1}, we get for all $u \in X$
\begin{eqnarray}
\Psi(\overline x_k) - \Psi(u)  
&\le&  (1- \gamma_k) [\Psi(\overline x_{k-1})  - \Psi(u)]
+ \gamma_k \left\{
\frac{\beta_k}{1 - P_{T_k}} [V(x_{k-1}, u) - V(x_k, u)] \right. \nonumber\\
 &&\quad \left. + \frac{P_{T_k}}{1 - P_{T_k} } 
\sum_{i=1}^{T_k} \frac{1}{p_i P_{i-1}} \left[\frac{\left(\tilde M + \|\delta_{k,i}\|_*\right)^2}{2\beta_k p_i } 
+ \langle \delta_{k,i}, u - u_{k,i-1} \rangle \right]\right\}.
\end{eqnarray}

Using the above inequality and Lemma~\ref{lem:lemma3.17_lan}, we conclude that for all $u \in X$
\begin{eqnarray}
\Psi(\overline x_N) - \Psi(u)   &\le& \Gamma_N (1- \gamma_1) [\Psi(\overline x_{0})  - \Psi(u)] \nonumber\\
&&\quad + \Gamma_N 
\sum_{k=1}^N \frac{\beta_k \gamma_k}{\Gamma_k (1 - P_{T_k})}
\left[
V(x_{k-1}, u) - V(x_k, u) \right] 
\nonumber \\
 &&\quad + \Gamma_N \sum_{k=1}^N 
 \frac{\gamma_k P_{T_k}}{\Gamma_k (1 - P_{T_k}) }\sum_{i=1}^{T_k} \frac{1}{p_i P_{i-1}} \left[\frac{\left(\tilde M + \|\delta_{k,i}\|_*\right)^2}{2\nu \beta_k p_i } 
+ \langle \delta_{k,i}, u - u_{k,i-1} \rangle \right].\label{eq:technical_sliding}
\end{eqnarray}

From \eqref{gamma_k_main} it follows that for all $u \in X$
\begin{eqnarray}
\sum_{k=1}^N \frac{\beta_k\gamma_k}{\Gamma_k (1 - P_{T_k})}
\left[
V(x_{k-1}, u) -  V(x_k, u) \right] &\le& \frac{\beta_1 \gamma_1  }{\Gamma_1 (1 - P_{T_1})} V(x_0, u) - \frac{\beta_N  \gamma_N }{\Gamma_N (1 - P_{T_N})} V(x_N, u) \notag\\
&\le& \frac{\beta_1}{1 - P_{T_1}} V(x_0, u), \label{bnd_dist1}
 \end{eqnarray}
where the last inequality follows from the facts that $\gamma_1 = \Gamma_1 = 1$, $P_{T_N} \le 1$, and $V(x_N, u) \ge0$. Inequality \eqref{bnd_dist1} and the fact that $\gamma_1 = 1$ together with inequality \eqref{eq:technical_sliding} imply that for all $u\in X$
    \begin{eqnarray}
    \label{1.69}
    \Psi(\overline x_N) - \Psi(u) &\leq& \frac{\beta_k }{1 - P_{T_1}} V(x_0,u)\notag\\
    &&\quad+
    \Gamma_N \sum\limits_{k=1}^N  \frac{\gamma_k P_{T_k}}{\Gamma_k(1 - P_{T_k})} \sum\limits_{i=1}^{T_k} \frac{1}{p_i P_{i-1}} \Bigg[ \frac{(\tilde M^2+ \|\delta_{k, i} \|^2_*)}{\beta_k p_i} + \langle \delta_{k,i}, u - u_{k,i-1}\rangle \Bigg].
    \end{eqnarray}
    
Next, we show that
\begin{equation}
    \label{temp3}  
    \mathbb{E}[||\delta_{k,i}||^2_{*}] \leq \sigma^2
\end{equation}
for $\sigma^2$ defined in \eqref{sig_main}. For all $x\in X$ we have
\begin{eqnarray}
    \label{sigg}
    \mathbb{E}[\|\delta\|^2_{*}] &=& \mathbb{E}[\|\tilde f_r'(x) - \nabla F(x)\|^2_{*}] \leq  2\mathbb{E}\|\tilde f_r'(x) \|^2_{*} + 2\mathbb{E}\|\nabla F(x)\|^2_{*} \nonumber\\ 
    &\overset{\eqref{exp_square_main}}{\leq}& 4p_*^2\left(cnM^2 + \frac{n^2\Delta^2 }{r^2} \right)+  2\|\nabla F(x)\|^2_{*} \nonumber\\
    &\overset{\eqref{eq:norm_nabla_F_bound_main}}{=}& 4p_*^2\left(cnM^2 + \frac{n^2\Delta^2 }{r^2} \right) + 2\tilde c^2 n p_*^2 M^2 \notag\\
    &=& 4p_*^2\left(CnM^2 + \frac{n^2\Delta^2 }{r^2} \right)\nonumber,
\end{eqnarray}
where $C \eqdef c + \nicefrac{\tilde{c}^2}{2}$. For the inner product we have the following bound:
\begin{equation}
    \label{temp2} 
    \mathbb{E}[\langle\delta_{k,i}, u - u_{k,i-1}\rangle] \overset{\eqref{tilde_f}}{=} \frac{n}{2r}\mathbb{E}[\langle\Delta_{k,i}e_{k,i}, u - u_{k,i-1}\rangle] \le \frac{n}{2r}\mathbb{E}[|\Delta_{k,i}|\cdot\|e_{k,i}\|_*\cdot \|u - u_{k,i-1}\|] \overset{\eqref{delta},\eqref{eq:condition_on_u}}{\le} \frac{\Delta n D_Xp_*}{r}.
\end{equation}
Taking mathematical expectation from the both sides of \eqref{1.69} and using \eqref{temp3} and \eqref{temp2} we obtain \eqref{t2_1}.

\subsection{Proof of Corollary~\ref{cor:main}}
    Using recurrences \eqref{p_t_main} and \eqref{pt_main} we obtain
    \begin{equation}
    \label{Pt1}
    P_t = \frac{2}{(t+1)(t+2)},
    \end{equation}
    \begin{equation}
    \label{Pt2}
    P_{T_k} \leq P_{T_{k-1}} \leq \ldots \leq P_{T_1} \leq \frac{1}{3}
    \end{equation}
    and from relations \eqref{gamma_kk_main} and \eqref{k_main} we derive that
    \begin{equation}
    \label{gk}
    \Gamma_k = \frac{2}{k(k+1)},
    \end{equation}
    which implies \eqref{gamma_main}. 
    
    From \eqref{k_main}, \eqref{Pt1}, \eqref{Pt2} we derive \eqref{gamma_k_main}. Simple calculations and relations \eqref{pt_main}, \eqref{Pt1} imply
    \begin{equation}
    \label{p1}
    \sum\limits_{i=1}^{T_k} \frac{1}{p_i^2 P_{i-1}} = 2 \sum\limits_{i=1}^{T_k} \frac{i+1}{i} \leq 4 T_k,
    \end{equation}
    \begin{equation}
    \label{p2}
    \sum\limits_{i=1}^{T_k} \frac{1}{p_i P_{i-1}} = \frac{1}{2}\sum\limits_{i=1}^{T_k} i = \frac{1}{4} T_k(T_k + 1).
    \end{equation}
    Next, one can see from \eqref{k_main}, \eqref{gk}, \eqref{p1}, \eqref{p2} that
    \begin{equation}
    \label{p3}
    \sum\limits_{i=1}^{T_k} \frac{\gamma_k P_{T_k}}{\Gamma_k \beta_k (1 - P_{T_k}) p_i^2 P_{i-1}} \leq \frac{4 \gamma_k P_{T_k} T_k}{\Gamma_k \beta_k (1 - P_{T_k})} = \frac{4k^2}{L(T_k + 3)},
    \end{equation}
    \begin{equation}
    \label{p4}
    \sum\limits_{i=1}^{T_k} \frac{\gamma_k P_{T_k}}{\Gamma_k  (1 - P_{T_k}) p_i P_{i-1}} \leq \frac{\gamma_k P_{T_k} T_k (T_k + 1)}{4 \Gamma_k (1 - P_{T_k})} = \frac{(T_k + 1) k}{2(T_k + 3)}.
    \end{equation}

    Finally, inequalities  \eqref{t2_1}, \eqref{Pt2}, \eqref{gk}, \eqref{p3}, \eqref{p4} imply
    \begin{eqnarray}
    \EE[\Psi(\overline x_N) - \Psi(x^*)] &\leq& \frac{2L}{N(N+1)}[3V(x_0, x^{*}) + 4\tilde D] + \frac{2n \Delta D_Xp_*}{r N(N+1)} \sum\limits_{k=1}^{N} \frac{(T_k + 1) k}{(T_k + 3)} \nonumber\\ 
    &\leq& \frac{2L}{N(N+1)}[3V(x_0, x^{*}) + 4\tilde D] + \frac{2n \Delta D_Xp_*}{r N(N+1)} \sum\limits_{k=1}^{N} k \nonumber\\
    &=& \frac{2L}{N(N+1)}[3V(x_0, x^{*}) + 4\tilde D] + \frac{n \Delta D_X p_*}{r}.\label{z1_main}
    \end{eqnarray} 

\subsection{Proof of Corollary~\ref{cor:main_corollary}}
The proof of \eqref{original_main} follows directly from \eqref{col1_main} and \eqref{orig_sm_main}. Using \eqref{original_main} and \eqref{r_delta_main} we get \eqref{bound_out_main}. Finally,
 \begin{eqnarray}
    \label{temp65}
    \sum\limits_{i=1}^{N} T_k &\overset{\eqref{k_main}}{\leq}&  \sum\limits_{i=1}^{N} \left(\frac{4N(\tilde M^2 + \sigma^2)k^2}{3D_{X,V}^2L^2} + 1\right) 
    \nonumber\\
    &=& \frac{2}{18} \frac{N^2(N+1)(2N+1)(\tilde M^2 + \sigma^2)}{D_{X,V}^2 L^2} + N \notag\\
    &\leq& \frac{2}{9} \frac{(N+1)^4(\tilde M^2 + \sigma^2)}{D_{X,V}^2 L^2} + N 
    \end{eqnarray} 
 
and
\begin{eqnarray}
    \label{temp651}
    \sigma^2 &=& O\left(p_*^2\left(CnM^2 + \frac{n^2\Delta^2 }{r^2}\right)\right) 
    \nonumber\\
    &=& O \left(p_*^2nM^2 \right),
\end{eqnarray}
\begin{eqnarray}
    \label{temp652}
    \tilde M^2 &=& O\left(nC_1^2M^2\right),
\end{eqnarray}
where we used $\varepsilon = O\left(\sqrt{n}MD_X\right)$. From \eqref{bound_out_main} we know that $N = O\left(\sqrt{\nicefrac{L D_{X,V}^2}{\varepsilon}} \right)$.  Together with \eqref{temp65}, \eqref{temp651},  \eqref{temp652} it gives \eqref{bound_in_main}. 

\section{Missing Proofs from Section~\ref{sec:main_res_str_cvx}}
\subsection{Proof of Theorem~\ref{thm:mzosa_main}}
    We prove this result by induction. From \eqref{z1_main} we have
    \begin{eqnarray*}
    \EE[\Psi(y_i) - \Psi(y^*)\mid y_{i-1}] \leq \frac{2L}{N_0(N_0+1)}\left(3V(y_{i-1}, y^{*}) + 4\tilde D\right)+ \frac{n \Delta D_X p_*}{r}.
    \end{eqnarray*}
    Here we use that $y_{i-1}$ is $x_0$ and $y_i$ is output for {\tt M-zoSA} after $i$-th iteration. Since $\Psi$ is a sum of convex and $\mu$-strongly convex function we have that $\Psi$ is $\mu$-strongly convex and
    \begin{equation*}
        \EE[\Psi(y_i) - \Psi(y^*)\mid y_{i-1}] \leq \frac{2L}{N_0(N_0+1)}\left(\frac{3}{\mu} \left(\Psi(y_{i-1}) - \Psi(y^*)\right) + 4\tilde D\right) + \frac{n \Delta D_X p_*}{r}.
    \end{equation*}
    Taking the full expectation from the both sides of previous inequality, using the induction hypothesis and the definition of $\tilde D$, we conclude that
    \begin{equation*}
        \EE[\Psi(y_i) - \Psi(y^*)] \leq \frac{2L}{N_0^2}\frac{5 \rho_0}{\mu 2^{i-1}} + \frac{2L}{N_0^2}\frac{6 n \Delta D_X}{\mu r} + \frac{n \Delta D_X}{r} \leq \frac{\rho_0}{2^i} + \frac{2n \Delta D_X p_*}{r},
    \end{equation*}
    where the last inequality follows from the definition of $N_0$.

\subsection{Proof of Corollary~\ref{cor:mzosa_main}}
    Combining inequalities \eqref{orig_sm_main} and \eqref{converg_str_conv} we get \eqref{original_strconv}. Next, if $r$ and $\Delta$ satisfy \eqref{r_delta_strconv}, then $2rM + \frac{2n\Delta D_X p_*}{r} = O(\e)$. Taking the total number of phases, i.e.\ restarts, $I=\lceil\log_2 \max\left[ 1, \nicefrac{\rho_0}{\varepsilon}\right]\rceil$ we get $\frac{\rho_0}{2^I} = O(\e)$ and, as a consequence, $\EE\left[\Psi_0(y_i) - \Psi_0(y^*)\right] = O(\e)$. Therefore, the total number of $\nabla g$ computations, which equals $N_0 I$, is bounded by \eqref{bound_out_stconv}.
    
    Now let us derive a bound on the total number of $\tilde f_r'$ computations. Without loss of generality, let us assume that $\rho_0 > \varepsilon$. Using the previous bound on $I$ and the definition of $T_k$, we get
\begin{eqnarray*}
    \sum_{i=1}^I \sum_{k=1}^{N_0} T_k &\le& \sum_{i=1}^I \sum_{k=1}^{N_0} \left( \frac{\mu N_0 (\tilde M^2+\sigma^2) k^2}{\rho_0 L^2} 2^{i} +1\right) \\
    &\leq& \sum_{i=1}^I \left[\frac{\mu N_0 (\tilde M^2+\sigma^2)}{3 \rho_0 L^2}(N_0+1)^3 2^{i}  + N_0\right]\\
    & \leq& \frac{\mu N_0 (N_0 + 1)^3 (\tilde M^2+\sigma^2)}{3 \rho_0 L^2} 2^{I+1} + N_0 I\\
    &\leq& \frac{4 \mu (N_0 + 1)^4 (\tilde M^2+\sigma^2)}{3 \varepsilon L^2} + N_0 I.
\end{eqnarray*}
This inequality, the definition of $N_0$ and inequalities \eqref{temp651} and \eqref{temp652} give the bound \eqref{bound_in_stconv} for the total number of $\tilde f_r'$ computations.

\end{document}